\documentclass[pdflatex,sn-mathphys-num]{sn-jnl}

\usepackage{graphicx}%
\usepackage{multirow}%
\usepackage{amsmath,amssymb,amsfonts}%
\usepackage{amsthm}%
\usepackage{mathrsfs}%
\usepackage[title]{appendix}%
\usepackage{xcolor}%
\usepackage{textcomp}%
\usepackage{manyfoot}%
\usepackage{algorithm}%
\usepackage{algorithmicx}%
\usepackage{algpseudocode}%
\usepackage{listings}%
\usepackage{verbatim}
\usepackage{mathrsfs}
\usepackage{bm}
\usepackage{color}
\usepackage{xcolor}
\usepackage{caption}
\usepackage{subcaption}
\usepackage{enumitem}
\usepackage{cases}
\usepackage{multirow}
\usepackage{tikz}
\usepackage{booktabs}%
\usepackage{hhline}
\usetikzlibrary{matrix}
\usetikzlibrary{arrows}
\usetikzlibrary{calc}
\usepackage{pgfplots}
\usepackage{pifont}
\usepackage{ragged2e}
\usepackage{array}
\usepackage{tabularx}

\theoremstyle{thmstyleone}%
\newtheorem{theorem}{Theorem}
\newtheorem{proposition}[theorem]{Proposition}%

\newtheorem{corollary}[theorem]{Corollary}%

\theoremstyle{thmstyletwo}%
\newtheorem{remark}{Remark}%

\theoremstyle{thmstylethree}%

\usepackage[colorinlistoftodos,prependcaption,backgroundcolor=black!5!white,bordercolor=red]{todonotes}

\renewcommand{\nu}{\vec{\mathbf{n}}}
\newcommand{\HH}{\mathcal{H}}
\newcommand{\scal}[2]{{\left\langle{{#1}\mid{#2}}\right\rangle}}
\newcommand{\menge}[2]{\big\{{#1}~\big |~{#2}\big\}} 
\newcommand{\Id}{\ensuremath{\operatorname{Id}}}
\newcommand{\zer}{\ensuremath{\operatorname{zer}}}
\newcommand{\dom}{\ensuremath{\operatorname{dom}}}

\newcommand{\prox}{{\rm prox}}
\newcommand{\gra}{{\rm gra\,}}

\newcommand{\RX}{]-\infty,+\infty]}
\newcommand{\RPP}{\ensuremath{\left]0,+\infty\right[}}
\newcommand{\RP}{\ensuremath{\left[0,+\infty\right[}}
\newcommand{\Fix}{\ensuremath{\text{\rm Fix}\,}}

\newcommand{\RR}{\mathbb{R}}
\newcommand{\NN}{\mathbb{N}}

\newtheorem{problem}[theorem]{Problem}
\graphicspath{
	{figures/},
	{figures/u-v/},
}
\makeatletter
	\def\input@path{%
		{figures/},%
	}
\makeatother

\renewcommand{\leq}{\leqslant}
\renewcommand{\geq}{\geqslant}
\renewcommand{\le}{\leqslant}
\renewcommand{\ge}{\geqslant}
\newcommand{\cA}{\mathcal{A}}
\newcommand{\cB}{\mathcal{B}}
\newcommand{\cM}{\mathcal{M}}
\newcommand{\rate}{\omega}

\DeclareFontEncoding{FMS}{}{}
\DeclareFontSubstitution{FMS}{futm}{m}{n}
\DeclareFontEncoding{FMX}{}{}
\DeclareFontSubstitution{FMX}{futm}{m}{n}
\DeclareSymbolFont{fouriersymbols}{FMS}{futm}{m}{n}
\DeclareSymbolFont{fourierlargesymbols}{FMX}{futm}{m}{n}
\DeclareMathDelimiter{\nr}{\mathord}{fouriersymbols}{152}{fourierlargesymbols}{147}

\DeclareMathOperator*{\argmin}{arg\,min}
\DeclareMathDelimiter{\nr}{\mathord}{fouriersymbols}{152}{fourierlargesymbols}{147}
\DeclareMathAlphabet{\mathpzc}{OT1}{pzc}{m}{it}

\begin{document}

\title[Douglas--Rachford algorithm with relaxed resolvents]{Douglas--Rachford algorithm with relaxed resolvents: convergence and tight linear rates}


\author[1]{\fnm{Luis} \sur{Briceño-Arias}}\email{luis.briceno@usm.cl}
\equalcont{These authors contributed equally to this work.}

\author*[2]{\fnm{Fernando} \sur{Roldán}}\email{fernandoroldan@udec.cl}
\equalcont{These authors contributed equally to this work.}

\author[3]{\fnm{Adrien} \sur{Taylor}}\email{adrien.taylor@inria.fr}
\equalcont{These authors contributed equally to this work.}

\affil[1]{\orgdiv{Departamento de Matemática}, \orgname{Universidad Técnica Federico Santa María}, \orgaddress{\street{Av. Vicuña Mackenna}, \city{San Joaquín}, \postcode{8940897}, \state{Santiago}, \country{Chile}}}

\affil*[2]{\orgdiv{Departamento de Ingenieria Matemática and CI$^2$MA}, \orgname{Universidad de Concepción}, \orgaddress{\street{Avda. Esteban Iturra}, \city{Concepción}, \postcode{4030000}, \state{Biobío}, \country{Chile}}}

\affil[3]{\orgdiv{Thoth team}, \orgname{Inria Centre at the University Grenoble Alpes}, \orgaddress{\street{655 Avenue de l'Europe}, \city{Montbonnot-Saint-Martin}, \postcode{38330}, \state{Auvergne-Rhône-Alpes}, \country{France}}}


\abstract{In this paper we study the convergence of a generalization of the Douglas--Rachford algorithm for finding a zero of
the sum of two maximally monotone operators, where reflected resolvents are extended to relaxed resolvents. This extension 
allows for different step sizes in the resolvents under cocoercivity.
We show that the resulting method 
converges weakly when its parameters lie in a range determined by the 
cocoercivity constants, which reduces to the classical algorithm when 
these constants vanish. Our main results concern the cocoercive case when one of the operators is in addition 
strongly monotone, in which the composition is a strict contraction. For 
both monotone inclusions and convex minimization, we compute the 
contraction factor and its optimal parameters in closed form, and show 
that the resulting rates are tight. Whenever the second operator is 
cocoercive, the optimal rates are strictly smaller than the corresponding 
rates of forward-backward splitting and of the classical 
Peaceman--Rachford and Douglas--Rachford algorithms and, in convex 
minimization, than those of FISTA.
All theoretical rates are corroborated by the numerical worst cases obtained
via performance estimation.}





\maketitle

	\section{Introduction}
	In this paper we aim at solving the following problem defined in a real Hilbert space $\HH$. We recall that an operator $\mathcal{M}\colon\HH\to 2^{\HH}$ is $\beta$-cocoercive
	for some $\beta\ge 0$ if
	\begin{equation}
		\label{e:defcoco}
		(\forall (x,u)\in\gra\mathcal{M})(\forall (y,v)\in\gra\mathcal{M})\quad 
		\scal{x-y}{u-v}\ge\beta
		\|u-v\|^2.
	\end{equation}
	
	Note that $\mathcal{M}$ is single-valued if $\beta>0$ and  
	we say that $\mathcal{M}$ is monotone if it is $0$-cocoercive.
	\begin{problem}
		\label{prob:main}
		Let $\beta_{\cA}\ge 0$, let $\beta_{\cB}\ge 0$, and let 
		$\mathcal{A}\colon\HH\to 
		2^{\HH}$ and $\mathcal{B}\colon\HH\to 2^{\HH}$ 
		be maximally monotone operators such that $\mathcal{A}$ is 
		$\beta_{\cA}$-cocoercive and $\mathcal{B}$ is $\beta_{\cB}$-cocoercive.  
		The 
		problem is to 
		\begin{equation}
			\label{e:probintro}
			\text{find}\quad x\in\HH\quad\text{such that}\quad 
			0\in\mathcal{A}x+\mathcal{B}x,
		\end{equation}
		under the assumption that its solution set $Z$ is nonempty.
	\end{problem}
	
	A classical method to solve Problem~\ref{prob:main} is the 
	\textit{Douglas--Rachford splitting} (DRS) \cite{Doug56,Lion79}, 
	which 
	iterates
	\begin{align}
		(\forall n\in\NN)\quad
		&\left\lfloor 
		\begin{array}{ll}
			x_{n}=J_{\mathcal{A}}z_n\\
			y_n=J_{\mathcal{B}}\big(2x_{n}-z_n\big)\\
			z_{n+1}=z_n+2\lambda(y_n-x_{n}),
			\label{e:DR}
		\end{array}
		\right.
	\end{align}
	where $z_0\in\HH$, $\lambda\in\RPP$, and 
	$J_{\mathcal{M}}=(\Id+\mathcal{M})^{-1}$ is the \textit{resolvent of 
		$\mathcal{M}\colon\HH\to 2^{\HH}$} and $\Id$ is the identity operator in $\HH$. 
	Note that
	\eqref{e:DR} can be written equivalently as
	\begin{equation}
		(\forall n\in\NN)\quad 
		z_{n+1}=Tz_n,\quad\text{where}\:\: 
		T=(1-\lambda)\Id+\lambda\,
		(2J_{\mathcal{B}}-\Id)\circ (2J_{\mathcal{A}}-\Id)
	\end{equation}
	and, for every operator $\mathcal{M}\colon\HH\to 2^{\HH}$, 
	$2J_{\mathcal{M}}-\Id$ is the \textit{reflected 
		resolvent of $\mathcal{M}$}. In the case when $\mathcal{M}$ is 
	maximally monotone, we have that $2J_{\mathcal{M}}-\Id$ is 
	nonexpansive 
	\cite[Corollary~23.11(ii)]{bauschke2011:Convex_analysis}. Hence, 
	in the context of \eqref{e:probintro}, 
	$(2J_{\mathcal{B}}-\Id)\circ (2J_{\mathcal{A}}-\Id)$ is also nonexpansive 
	and, thus, when 
	$\lambda\in\left]0,1\right[$, $T$ is 
	$\lambda$-averaged nonexpansive 
	\cite[Definition~4.33]{bauschke2011:Convex_analysis} and we 
	obtain the weak convergence of DRS to a solution to 
	\eqref{e:probintro} 
	\cite[Theorem~26.11]{bauschke2011:Convex_analysis}.
	
	In the limit case when $\lambda= 1$, $T$ is merely 
	nonexpansive and the convergence of DRS is not guaranteed in 
	general. This case is known as the
	\textit{Peaceman--Rachford splitting} (PRS), which converges only 
	under additional assumptions such as the cocoercivity or the strong 
	monotonicity of 
	$\mathcal{B}$ \cite[Remark~2]{Lion79}. 
    When $\mathcal{A}$ or $\mathcal{B}$ is cocoercive, it can also be 
activated through an explicit step, which leads to the forward-backward 
splitting (FBS), and in smooth convex optimization to the proximal 
gradient method and its accelerated version FISTA \cite{BeckTeboulle2009}.
	A theoretical and numerical comparison between an explicit Euler 
scheme, FBS, PRS, and DRS when 
$\mathcal{A}$ and $\mathcal{B}$ are cocoercive and $\mathcal{A}$ is 
strongly monotone can be found in \cite{Sigpro}. In the case of 
strongly convex optimization problems, i.e., when $\mathcal{A}$ and 
$\mathcal{B}$ are gradients of convex differentiable functions $f$ 
and $g$ from $\HH$ to $\RR$ and $f$ is strongly convex, PRS has the 
smallest linear convergence rate in most cases, except when the 
cocoercivity constant of $\nabla g$ is large compared with that of 
$\nabla f$ \cite[Proposition~7]{Sigpro}. A direct comparison of the 
rates in \cite[Proposition~3]{Sigpro} and 
\cite[Theorem~5.6]{Giselsson2017} shows that the same holds in the 
monotone case, with $\nabla f$ and $\nabla g$ replaced by $\mathcal{A}$ 
and $\mathcal{B}$.
Hence, none of these methods has the best rate for all values of the 
constants, and the best choice depends on the problem. This raises 
the question of whether a single splitting method can improve all 
these rates by exploiting the cocoercivity of both operators.
 
	In this paper we provide an extension of the DRS for solving 
	\eqref{e:probintro} when $\mathcal{A}$ and $\mathcal{B}$ are 
	$\beta_{\cA}$- and $\beta_{\cB}$-cocoercive, respectively, for some $\beta_{\cA}\ge 0$ and $\beta_{\cB}\ge 0$. Our extension is based on the \textit{relaxed resolvent of 
		$\mathcal{M}$}, for some $\theta\ge-1$, defined by 
	\begin{equation}
		\label{e:relaxedrefl}
		R_\mathcal{M}^{\theta}=(1+\theta)J_\mathcal{M}-\theta\Id,
	\end{equation}
	which is a conic average of the resolvent in the sense of
	\cite{BauschkeMoursiWang2021}. Note that $R_\mathcal{M}^1=2J_\mathcal{M}-\Id$ 
	is the reflected resolvent of $\mathcal{M}$, $R_\mathcal{M}^{-1}=\Id$, and 
	$R_\mathcal{M}^0=J_\mathcal{M}$.
	
	Modifications of the reflected resolvent have been considered before. In 
	\cite{Arta19}, an alternating modified reflections algorithm is derived in 
	order to compute the resolvent of the sum of two maximally monotone 
	operators. Closer to our setting, relaxed resolvents of the form 
	\eqref{e:relaxedrefl} are the building block of the adaptive 
	Douglas--Rachford algorithm of \cite{DaoPhan2019}, whose iteration 
	coincides, up to a reparametrization, with the composition 
	$R_{\theta\mathcal{B}}^{1/\theta}\circ R_{\mathcal{A}}^{\theta}$ studied in this paper. The admissible parameters in \cite{DaoPhan2019} are governed by the strong and weak 
	monotonicity constants of $\mathcal{A}$ and $\mathcal{B}$, and linear rates 
	are obtained when $\mathcal{A}$ is Lipschitz continuous, or when 
	$\mathcal{B}$ is linear and Lipschitz continuous. On the other hand, 
	this composition has 
	recently been studied in 
	\cite{AlcantaraDaoTakeda2025} for maximally comonotone 
	operators. This notion reduces to cocoercivity when the corresponding 
	constant is nonnegative, via the Attouch--Th\'era duality and conical 
	averagedness \cite{BauschkeMoursiWang2021}, yielding weak convergence for 
	constant relaxations together with extensions to multioperator inclusions. 
	In \cite{NilssonAkermanGiselsson2025}, the relaxations for which the scheme 
	converges unconditionally are characterized for subdifferentials of 
	functions in $\Gamma_0(\HH)$ and without any cocoercivity assumption. In 
	particular, for merely maximally monotone operators, unconditional 
	convergence forces $\theta=1$. A different axis of variation is explored in 
	\cite{AtenasBauschkeDaoTam2026}, where the resolvent parameter is allowed to 
	vary along the iterations, whereas here the two operators are evaluated with 
	different but iteration-independent step-sizes.
	
	Our purpose is to study the role of cocoercivity in this scheme. We prove 
	that $R_\mathcal{M}^{\theta}$ is averaged nonexpansive when $\mathcal{M}$ is 
	$\beta$-cocoercive and $\theta\in\left]-1,1+2\beta\right]$, and that the 
	composition $R_{\theta\mathcal{B}}^{1/\theta}\circ R_{\mathcal{A}}^{\theta}$ 
	is averaged nonexpansive on a range of $\theta$ determined by $\beta_{\cA}$ and 
	$\beta_{\cB}$, which reduces to $\{1\}$ when $\beta_{\cA}=\beta_{\cB}=0$, in which case we 
	recover the classical DRS. Since the solutions to \eqref{e:probintro} are 
	obtained from the fixed points of this composition, the convergence of our 
	algorithm follows from properties of averaged nonexpansive mappings, and the 
	cocoercivity allows for relaxation steps beyond those available in the 
	classical case. Table~\ref{tab:relaxedresolvent} summarizes known and new properties of 
the relaxed resolvent, which serve as building blocks for the study of 
the composition.
	
\begin{figure}[htbp]
\centering
\includegraphics[width=\textwidth]{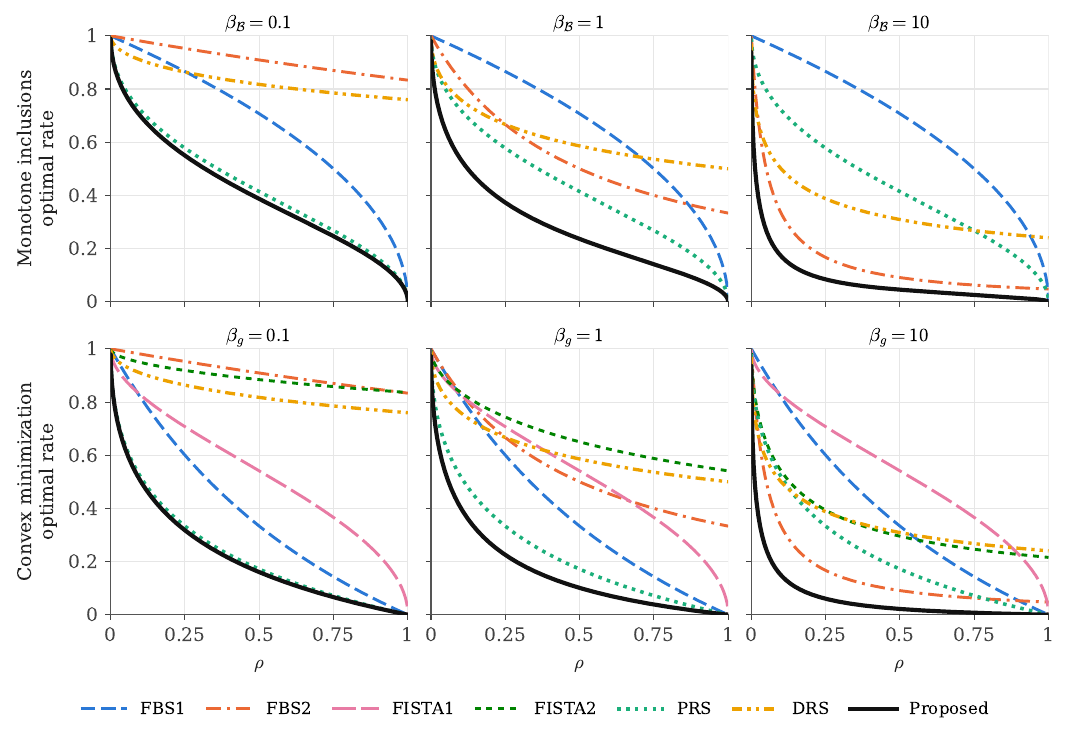}
\caption{Optimal linear rates in Table~\ref{tab:rates_comparison} as 
functions of $\rho\in\left]0,1\right]$, for $\beta_{\cA}=1$ (top) and 
$\beta_f=1$ (bottom), and for $\beta_{\cB}$ (top) and $\beta_g$ (bottom) 
in $\{0.1,1,10\}$.}
\label{fig:rates_comparison}
\end{figure}

	When $\mathcal{A}$ is in addition strongly monotone, the 
	composition is a strict contraction. We obtain a closed form of the contraction factor depending on the step-size and relaxation parameter, we minimize it over both constants, and we exhibit operators 
	attaining it, so that the resulting linear rates cannot be improved. By 
	contrast, the linear rates in \cite{DaoPhan2019} are upper bounds for a 
	relaxed iteration under Lipschitz-type assumptions that are not shown to be 
	tight, sharpness being established there only for the classical 
	Douglas--Rachford operator in a special class 
	\cite[Remark~4.11]{DaoPhan2019}. We obtain analogous results, with smaller constants, 
	for the minimization of the sum of two convex functions, where 
	the optimal rate is strictly smaller than the best rate available for the 
	classical Douglas--Rachford algorithm in the same setting. All theoretical 
	rates are corroborated by the numerical worst cases obtained via 
	performance estimation. In Table~\ref{tab:composition} we summarize the results in the literature and our contributions on linear rates for fixed-point iterations governed by $R_{\theta\tau\cB}^{1/\theta}\circ R_{\tau\cA}^{\theta}$. Table~\ref{tab:rates_comparison} compares the optimal rates obtained 
in this paper with those of FBS, PRS, and DRS in the strongly 
monotone and cocoercive setting and also FISTA for strongly convex smooth minimization problems. Our optimal rates are strictly smaller 
than those of PRS and DRS whenever $\beta_{\cB}>0$, and 
Figure~\ref{fig:rates_comparison} shows that they are also smaller 
than the optimal rates of FBS and the best known rates of FISTA in all 
the configurations considered, at the price of two resolvent 
evaluations per iteration.

\begin{table}[!htb]
\centering
\caption{\scriptsize Properties of $R^{\theta}_{\cM}$ for a maximally monotone 
operator $\cM\colon\HH\to2^{\HH}$. Sharp: \ding{52} the constant is attained 
by an operator of the class. $^\dagger$ Can be deduced from the arguments of this 
paper.}
\label{tab:relaxedresolvent}
{\scriptsize\renewcommand{\arraystretch}{1.3}\setlength{\tabcolsep}{3pt}
\begin{tabular}{@{}
  >{\raggedright\arraybackslash}p{0.19\textwidth}
  >{\raggedright\arraybackslash}p{0.12\textwidth}
  >{\raggedright\arraybackslash}p{0.24\textwidth}
  >{\raggedright\arraybackslash}p{0.13\textwidth}
  >{\raggedright\arraybackslash}p{0.25\textwidth}@{}}
\specialrule{2pt}{1pt}{1pt}
\bfseries Assumption on $\cM$ & \bfseries Range of $\theta$ 
& \bfseries Property of $R^{\theta}_{\cM}$ & \bfseries Sharp 
& \bfseries Reference\\
\hline
maximally monotone & $\left]-1,1\right]$ 
& $\frac{1+\theta}{2}$-averaged nonexpansive 
& \ding{52}$^\dagger$ ($\cM=N_{\{0\}}$) 
& \cite[Corollary~23.11(ii)]{bauschke2011:Convex_analysis} ($\theta=1$), 
Prop.~\ref{t:genDRtau}\ref{t:genDRtau2+} with $\beta_{\cA}=0$\\
 & $\left]1,+\infty\right[$ 
& $\theta$-Lipschitz, not nonexpansive 
& \ding{52}$^\dagger$ ($\cM=N_{\{0\}}$) 
& $^\dagger$ Case~3 in the proof of Prop.~\ref{p:reflectcontr}\\
$\beta$-cocoercive, $\beta>0$ & $\left]-1,1+2\beta\right]$ 
& $\frac{1+\theta}{2(1+\beta)}$-averaged nonexpansive 
& \ding{52}$^\dagger$ ($\cM=\beta^{-1}\Id$) 
& Prop.~\ref{t:genDRtau}\ref{t:genDRtau2+}, 
\cite[Proposition~5.3]{Giselsson2017} ($\theta=1$)\\
 & $\left]1+2\beta,+\infty\right[$ 
& $\frac{\theta-\beta}{1+\beta}$-Lipschitz, not nonexpansive 
& \ding{52} $^\dagger$ ($\cM=\beta^{-1}\Id$) 
& Thm.~\ref{t:rrcontractive}\ref{t:rrcontractive1}{\&}\ref{t:rrcontractive3}\\
$\rho$-strongly monotone and $\beta$-cocoercive, $\rho>0$, $\beta\ge0$, 
$\beta\rho\le1$ 
& $\left[-1,+\infty\right[$ 
& $\rate_{\cM}(\theta)$-Lipschitz, see \eqref{e:defrhoa}, contraction if 
and only if $\theta\in\left]-1,1+2\beta\right[$ 
& \ding{52} ($\beta>0$) 
& Thm.~\ref{t:rrcontractive}\ref{t:rrcontractive1}{\&}\ref{t:rrcontractive3} with $\tau=1$, Prop.~\ref{p:reflectcontr}, 
cf.\ \cite[Theorem~7.2]{Giselsson2017}\\
$\partial f$, $f$ $\rho$-strongly convex, $\nabla f$ $\beta$-cocoercive 
& $\left[-1,+\infty\right[$ 
& $\max\big\{\frac{1-\theta\rho}{1+\rho},\frac{\theta-\beta}{1+\beta}\big\}$-Lipschitz 
& \ding{52}$^\dagger$ (quadratic $f$) 
& Prop.~\ref{p:optiref}\ref{p:optirefi}, 
\cite[Theorem~1]{Giselsson2017IEEE} for 
$\theta\in\left]-1,1+2\beta\right[$\\
\specialrule{2pt}{1pt}{1pt}
\end{tabular}}
\end{table}
\begin{table}[!htbp]
\centering
\caption{\scriptsize Properties of 
$T_{\theta,\tau}=R_{\theta\tau\cB}^{1/\theta}\circ R_{\tau\cA}^{\theta}$, 
$\theta>0$, $\tau>0$, and of the associated relaxed iterations. For 
$\theta=1$, $T_{1,\tau}=R_{\tau\cB}\circ R_{\tau\cA}$ (Peaceman--Rachford). 
Sharp: \ding{52} the constant or the condition cannot be improved over the 
class, -- not addressed.}
\label{tab:composition}
{\scriptsize\renewcommand{\arraystretch}{1.3}\setlength{\tabcolsep}{3pt}
\begin{tabular}{@{}
  >{\raggedright\arraybackslash}p{0.15\textwidth}
  >{\raggedright\arraybackslash}p{0.13\textwidth}
  >{\raggedright\arraybackslash}p{0.16\textwidth}
  >{\raggedright\arraybackslash}p{0.22\textwidth}
  >{\raggedright\arraybackslash}p{0.08\textwidth}
  >{\raggedright\arraybackslash}p{0.17\textwidth}@{}}
\specialrule{2pt}{1pt}{1pt}
\bfseries $\cA$ & \bfseries $\cB$ 
& \bfseries Parameters & \bfseries Property & \bfseries Sharp 
& \bfseries Reference\\
\hline
\multicolumn{6}{@{}l}{\itshape Linear convergence}\\
$\rho$-strongly monotone, $\beta_{\cA}$-cocoercive 
& $\beta_{\cB}$-cocoercive 
& $\theta\in I(\tau)$, see \eqref{e:I(tau)} 
& $\rate(\theta,\tau)$-contraction \eqref{e:defomega}, optimal 
parameters \eqref{e:main}--\eqref{e:argmin} 
& \ding{52} ($\cB=0$) 
& Thm.~\ref{t:rrcontractive}\ref{t:rrcontractive2} and 
Thm.~\ref{t:rrcontractive}\ref{t:rrcontractive3}, 
Thm.~\ref{t:strongconvergenceDR}\\
$\rho$-strongly monotone, $\beta_{\cA}$-cocoercive 
& maximally monotone 
& $\theta=1$ 
& $\rate(1,\tau)$-contraction, optimal $\tau$ 
& \ding{52} 
& \cite[Theorem~7.4]{Giselsson2017}\\
$\rho$-strongly monotone & $\beta_{\cB}$-cocoercive 
& $\theta=1$, relaxed 
& rate $1/(1+\sqrt{\beta_{\cB}\rho})$ for a specific relaxation and step-size 
& \ding{52} 
& \cite[Theorem~5.6]{Giselsson2017}\\
$\rho$-strongly monotone, Lipschitz & maximally monotone 
& $\theta\in[(1+2\tau\rho)^{-1},1]$, relaxation $\kappa\in\left]0,1\right[$ 
& linear convergence 
& --
& \cite[Theorem~4.8, Corollary~4.9]{DaoPhan2019}\\
\multicolumn{2}{@{}>{\raggedright\arraybackslash}p{0.28\textwidth}}{combinations of strong monotonicity, 
cocoercivity, and Lipschitz continuity} 
& $\theta=1$, relaxed 
& tight contraction factors via performance estimation 
& \ding{52} 
& \cite{RyuTaylorBergelingGiselsson2020}\\
Lipschitz & strongly monotone 
& $\theta=1$ 
& linear convergence 
& -- 
& \cite{MoursiVandenberghe2019}\\
\hline
\multicolumn{6}{@{}l}{\itshape Convex optimization ($\cA=\partial f$, $\cB=\partial g$)}\\
$f$ $\rho$-strongly convex, $\nabla f$ $\beta_f$-cocoercive 
& $\nabla g$ $\beta_g$-cocoercive 
& $(1-\theta)\tau\in]-2\beta_f,2\beta_g]$ 
& $\widehat{\rate}(\theta,\tau)$-contraction \eqref{e:r1}, optimal value 
\eqref{e:main} 
& \ding{52} ($\beta_f\beta_g>0$) 
& Thm.~\ref{p:optiref}\ref{p:optirefii}--\ref{p:optirefiii}\\
$f$ $\rho$-strongly convex, $\nabla f$ $\beta_f$-cocoercive 
& $g\in\Gamma_0(\HH)$ 
& $\theta=1$ 
& $\widehat{\rate}(1,\tau)$-contraction \eqref{e:r1}, optimal value 
\eqref{e:main} with $\beta_g=0$ 
& \ding{52}
& \cite[Theorem~1, Proposition~3]{Giselsson2017IEEE}, 
\cite[Theorem~3.4(3)]{BricenoRoldan2026}\\
\specialrule{2pt}{1pt}{1pt}
\end{tabular}}
\end{table}

\begin{table}[!htbp]
\centering
\caption{\scriptsize Optimal linear rates in the strongly monotone and 
cocoercive setting. Monotone inclusions: $\cA$ is $\rho$-strongly monotone 
and $\beta_{\cA}$-cocoercive, and $\cB$ is $\beta_{\cB}$-cocoercive. Convex 
minimization of $f+g$: $f$ is $\rho$-strongly convex, $\nabla f$ is 
$\beta_f^{-1}$-Lipschitz continuous, and $\nabla g$ is 
$\beta_g^{-1}$-Lipschitz continuous. FBS1 and FISTA1 activate explicitly 
the strongly monotone operator or the strongly convex function, and FBS2 
and FISTA2 the other one. In the proposed method the rates are attained with the parameters given in Theorem~\ref{t:rrcontractive} and Proposition~\ref{t:opti}. 
Sharp: \ding{52} the optimal rate is attained in the class. 
$^\dagger$ Follows from scalar or planar examples as in this paper. 
$^{*}$ Square root of the rate of a Lyapunov function that bounds 
$\|x_n-\widehat{x}\|^2$, hence an R-linear rate for $\|x_n-\widehat{x}\|$.}
\label{tab:rates_comparison}
{\scriptsize\renewcommand{\arraystretch}{1.5}\setlength{\tabcolsep}{3pt}
\begin{tabular}{@{}
  >{\raggedright\arraybackslash}p{0.27\textwidth}
  >{\raggedright\arraybackslash}p{0.33\textwidth}
  >{\raggedright\arraybackslash}p{0.07\textwidth}
  >{\raggedright\arraybackslash}p{0.27\textwidth}@{}}
\specialrule{2pt}{1pt}{1pt}
\bfseries Algorithm & \bfseries Optimal rate 
& \bfseries Sharp & \bfseries References\\
\hline
\multicolumn{4}{@{}l}{\itshape Monotone inclusions}\\
FBS1: $J_{\tau\cB}\circ(\Id-\tau\cA)$ 
& $\sqrt{1-\beta_{\cA}\rho}$ 
& \ding{52}$^\dagger$ 
& \cite{Sigpro}\\\addlinespace[4pt]
FBS2: $J_{\tau\cA}\circ(\Id-\tau\cB)$ 
& $\dfrac{1}{1+2\beta_{\cB}\rho}$ 
& \ding{52}$^\dagger$ 
& \cite{Sigpro}\\\addlinespace[4pt]
PRS 
& $\sqrt{\dfrac{1-\sqrt{\beta_{\cA}\rho}}{1+\sqrt{\beta_{\cA}\rho}}}$ 
& \ding{52} 
& \cite[Theorem~7.4]{Giselsson2017}, 
Thm.~\ref{t:rrcontractive}\\\addlinespace[4pt]
DRS 
& $\dfrac{1}{1+\sqrt{\beta_{\cB}\rho}}$ 
& \ding{52} 
& \cite[Theorem~5.6]{Giselsson2017}\\\addlinespace[4pt]
Proposed with \\[-.2cm]
$\widehat{\beta}=\max\{\beta_{\cA},\beta_{\cB}/(1+2\beta_{\cB}\rho)\}$\\[-.3cm]
$\widehat{\tau}=\beta_{\cB}+\sqrt{\beta_{\cB}^2+2\widehat{\beta}\beta_{\cB}
+\widehat{\beta}/\rho}$
& $\dfrac{\sqrt{\widehat{\beta}(1/\rho-\widehat{\beta})}}
{\widehat{\tau}+\widehat{\beta}}$ 
& \ding{52} 
& Thm.~\ref{t:rrcontractive}\\\addlinespace[4pt]
\hline
\multicolumn{4}{@{}l}{\itshape Convex minimization}\\
FBS1: $\prox_{\tau g}\circ(\Id-\tau\nabla f)$ 
& $\dfrac{1-\beta_f\rho}{1+\beta_f\rho}$ 
& \ding{52} 
& \cite{Sigpro,TaylorHendrickxGlineur2018}\\\addlinespace[4pt]
FBS2: $\prox_{\tau f}\circ(\Id-\tau\nabla g)$ 
& $\dfrac{1}{1+2\beta_g\rho}$ 
& \ding{52}$^\dagger$ 
& \cite{Sigpro,BricenoFISTA}\\\addlinespace[4pt]
FISTA1: $\nabla f$ and $\prox_{g}$ 
& $\sqrt{1-\sqrt{\beta_f\rho}}\;{}^{*}$ 
& -- 
& \cite[Theorem~10.42]{Beck2017}, \cite{BricenoFISTA}\\\addlinespace[4pt]
FISTA2: $\nabla g$ and $\prox_{f}$ 
& $\sqrt{1-\sqrt{\dfrac{\beta_g\rho}{1+\beta_g\rho}}}\;{}^{*}$ 
& -- 
& \cite[Theorem~4.10]{ChambollePock2016}, \cite{BricenoFISTA}\\\addlinespace[4pt]
PRS 
& $\dfrac{1-\sqrt{\beta_f\rho}}{1+\sqrt{\beta_f\rho}}$ 
& \ding{52} 
& \cite[Theorem~1, Proposition~3]{Giselsson2017IEEE}\\\addlinespace[4pt]
DRS 
& $\dfrac{1}{1+\sqrt{\beta_g\rho}}$ 
& \ding{52} 
& \cite[Theorem~5.6]{Giselsson2017}\\\addlinespace[4pt]
Proposed 
& $\dfrac{\sqrt{1+\beta_g\rho}-\sqrt{\rho(\beta_f+\beta_g)}}
{\sqrt{1+\beta_g\rho}+\sqrt{\rho(\beta_f+\beta_g)}}$ 
& \ding{52} 
& Prop.~\ref{t:opti}\\\addlinespace[4pt]
\specialrule{2pt}{1pt}{1pt}
\end{tabular}}
\end{table}

	The paper is organized as follows. In Section~\ref{sec:prel}
	we provide the preliminaries and the main properties we use throughout the paper. We study the convergence of our algorithm for monotone inclusions in Section~\ref{sec:moninclu} and for convex optimization in Section~\ref{sec:convex}, including linear rates under strong monotonicity and strong convexity, respectively.

	\section{Preliminaries}
	\label{sec:prel}
	Throughout this paper, $\HH$ is a real Hilbert space with inner product 
	$\scal{\cdot}{\cdot}$ and associated norm $\|\cdot\|=\sqrt{\scal{\cdot}{\cdot}}$. 
	We denote by $\Id$ the identity operator on $\HH$, and by $\rightharpoonup$ and 
	$\to$ the weak and the strong convergence in $\HH$, respectively. 
	
	Let $T\colon\HH\to\HH$. The set of fixed points of $T$ is 
	$\Fix T=\menge{x\in\HH}{Tx=x}$. Given $L\ge 0$, $T$ is $L$-Lipschitz continuous if, for every $x\in\HH$ and $y\in\HH$,
	$\|Tx-Ty\|\le L\|x-y\|,$
	and it is a strict contraction if it is $L$-Lipschitz continuous for some 
	$L\in\left[0,1\right[$.
	The operator $T$ is $\alpha$-averaged nonexpansive
	for some $\alpha\in\left]0,1\right]$ if 
	\begin{equation}
		\label{e:averaged}
		(\forall x\in\HH)(\forall y\in\HH)\quad 
		\|Tx-Ty\|^2\le\|x-y\|^2-\left(\frac{1-\alpha}{\alpha}\right)\|(\Id-T)x-(\Id-T)y\|^2.
	\end{equation} 
	Equivalently, $T$ is $\alpha$-averaged nonexpansive if and only if there exists 
	a nonexpansive operator $N\colon\HH\to\HH$ such that $T=(1-\alpha)\Id+\alpha N$ 
	\cite[Proposition~4.35]{bauschke2011:Convex_analysis}. Moreover, $T$ is nonexpansive if it is $1$-averaged nonexpansive, that is, 
	$1$-Lipschitz continuous, and it is firmly nonexpansive if it is 
	$\frac{1}{2}$-averaged nonexpansive.
	
	Let $\mathcal{A}\colon\HH\to 2^{\HH}$ be a set-valued operator. The domain, 
	the graph, and the set of zeros of $\mathcal{A}$ are, respectively,
	$\dom\mathcal{A}=\menge{x\in\HH}{\mathcal{A}x\neq\varnothing}$, 
	$\zer\mathcal{A}=\menge{x\in\HH}{0\in\mathcal{A}x}$, and $\gra\mathcal{A}=\menge{(x,u)\in\HH\times\HH}{u\in\mathcal{A}x}$.
	The inverse of $\mathcal{A}$ is the operator 
	$\mathcal{A}^{-1}\colon\HH\to 2^{\HH}$ with graph 
	$\menge{(u,x)\in\HH\times\HH}{(x,u)\in\gra\mathcal{A}}$, and the resolvent of 
	$\mathcal{A}$ is $J_{\mathcal{A}}=(\Id+\mathcal{A})^{-1}$. Let $\rho\ge 0$. 
	The operator $\mathcal{A}$ is $\rho$-strongly monotone if
	\begin{equation}
		\label{e:defstrong}
		(\forall (x,u)\in\gra\mathcal{A})(\forall (y,v)\in\gra\mathcal{A})\quad
		\scal{x-y}{u-v}\ge\rho\|x-y\|^2.
	\end{equation}
	The operator $\mathcal{A}$ is monotone if and only if it is $0$-cocoercive 
	(see \eqref{e:defcoco}), or, equivalently, $0$-strongly monotone, and it is 
	maximally monotone if it is monotone and, for every monotone 
	operator $\mathcal{B}\colon\HH\to 2^{\HH}$ such that 
	$\gra\mathcal{A}\subset \gra\mathcal{B}$, we have $\gra\mathcal{A}= 
	\gra\mathcal{B}$. If $\mathcal{A}$ is maximally monotone, then, for every 
	$\gamma>0$, $J_{\gamma\mathcal{A}}\colon\HH\to\HH$ is single-valued, firmly 
	nonexpansive, and $\Fix J_{\gamma\mathcal{A}}=\zer\mathcal{A}$, and the reflected 
	resolvent $R_{\gamma\mathcal{A}}=2J_{\gamma\mathcal{A}}-\Id$ is nonexpansive 
	\cite[Corollary~23.11]{bauschke2011:Convex_analysis}.
	
	Let $\mathcal{A}\colon\HH\to 2^{\HH}$ be maximally monotone. It follows from 
	\cite[Proposition~3.11(ii)]{BauschkeMoursiWang2021} that, for every 
	$\beta_{\mathcal{A}}\ge 0$,
	\begin{equation}
		\label{e:equivcoco}
		\text{$\mathcal A$ is $\beta_{\mathcal{A}}$-cocoercive
			\:\: $\Leftrightarrow$ \:\: $J_{\mathcal{A}}$ is $\frac{1}{2(1+\beta_{\mathcal{A}})}$-averaged nonexpansive},
	\end{equation}
	and from \cite[Proposition~23.13]{bauschke2011:Convex_analysis} that, for every 
	$\rho\ge 0$,
	\begin{equation}
		\label{e:equivstrong}
		\text{$\mathcal A$ is $\rho$-strongly monotone
			\quad $\Leftrightarrow$ \quad $J_{\mathcal{A}}$ is $(1+\rho)$-cocoercive}.
	\end{equation}
	
	Let $f\colon\HH\to\RX$. The domain 
	of $f$ is $\dom f=\menge{x\in\HH}{f(x)<+\infty}$, and $f$ is proper if 
	$\dom f\neq\varnothing$. We denote by $\Gamma_0(\HH)$ the class of all proper, 
	lower semicontinuous, convex functions from $\HH$ to $\RX$. Let 
	$f\in\Gamma_0(\HH)$. The Fenchel conjugate of $f$ is the function 
	$f^*\in\Gamma_0(\HH)$ defined by 
	$f^*\colon u\mapsto \sup_{x\in\HH} \left(\scal{x}{u}-f(x)\right)$.
	The subdifferential of $f$ is the operator $\partial f\colon x\mapsto \menge{u\in\HH}{(\forall y\in\HH)\:\: f(x)+\scal{y-x}{u}\le f(y)}$,
	it is maximally monotone, it satisfies 
	$(\partial f)^{-1}=\partial f^*$, and its set of zeros $\zer\partial f$ 
	coincides with the set of global minimizers of $f$, denoted by 
	$\argmin_{x\in \HH}f$. If $f$ is Fr\'echet differentiable at $x\in\HH$, then 
	$\partial f(x)=\{\nabla f(x)\}$. Given $\rho\ge 0$, $f$ is $\rho$-strongly 
	convex if $f-\frac{\rho}{2}\|\cdot\|^2$ is convex, in which case $\partial f$ is 
	$\rho$-strongly monotone. Moreover, if $f\colon\HH\to\mathbb{R}$ is convex and 
	Fr\'echet differentiable and $\beta>0$, the Baillon--Haddad theorem asserts that 
	$\nabla f$ is $\beta$-cocoercive if and only if $\nabla f$ is 
	$\beta^{-1}$-Lipschitz continuous 
	\cite[Corollary~18.17]{bauschke2011:Convex_analysis}.
	The proximity operator of $f$ is the single-valued mapping
	\begin{equation*}
		\prox_{f}\colon x\mapsto\argmin_{y\in\HH}\left(f(y)+\frac{1}{2}\|x-y\|^2\right),
	\end{equation*}
	which coincides with the resolvent $J_{\partial f}$. Hence, for every 
	$\gamma>0$, $\prox_{\gamma f}=J_{\gamma\partial f}$ is firmly nonexpansive and 
	$\Fix\prox_{\gamma f}=\argmin_{x\in\HH}f$. Given a nonempty closed convex 
	set $C\subset\HH$, we denote by $\iota_C$ its indicator function, which takes 
	the value $0$ on $C$ and $+\infty$ otherwise, by $N_C=\partial\iota_C$ its 
	normal cone operator, and by $P_C=\prox_{\iota_C}=J_{N_C}$ the projector onto 
	$C$. For instance, $\gra N_{\{0\}}=\{0\}\times\HH$ and $J_{N_{\{0\}}}\colon x\mapsto 0$.
	Finally, for every $f$ and $g$ in $\Gamma_0(\HH)$, we have 
	$\zer(\partial f+\partial g)\subset\argmin_{x\in\HH}(f+g)$, and equality holds, 
	e.g., if $\dom g=\HH$ \cite[Corollary~16.48]{bauschke2011:Convex_analysis}. 
	For a comprehensive treatment of monotone operator theory and convex analysis, 
	we refer the reader to \cite{bauschke2011:Convex_analysis}.

	\section{Monotone inclusions}\label{sec:moninclu}
	In this section we first show how the solutions to Problem~\ref{prob:main} are recovered from the fixed points of the composition of relaxed resolvents. Next, we study the averagedness and contraction properties of this composition.
	
	\begin{proposition}
		\label{t:genDRtau}
		In the context of Problem~\ref{prob:main}, the 
		following hold:
		\begin{enumerate}
			\item\label{t:genDRtau1} For every $\theta>0$, 
			$Z=J_{\mathcal{A}}(\Fix (R_{\theta\mathcal{B}}^{1/\theta}\circ 
			R_{\mathcal{A}}^{\theta}))$.
			
			\item\label{t:genDRtau2}
			For every $\theta\ge -1$, we have, for every 
			$x\in\HH$ and $y\in\HH$,
			\begin{equation}
				\label{e:ineqrelref}
				\|R_\mathcal{A}^{\theta}x-R_\mathcal{A}^{\theta}y\|^2
				\le 
				\|x-y\|^2-(1+\theta)(1+2\beta_{\cA}-\theta)
				\|(\Id-J_{\mathcal{A}})x-(\Id-J_\mathcal{A})y\|^2.
			\end{equation}
			
			\item\label{t:genDRtau2+} 
			For every $\theta\in\left]-1,1+2\beta_{\cA}\right]$, 
			$R_\mathcal{A}^{\theta}$ is $\alpha_{\mathcal{A}}(\theta)$-averaged 
			nonexpansive, where 
			\begin{equation}
				\alpha_{\mathcal{A}}(\theta)=\dfrac{1+\theta}{2(1+\beta_{\cA})}\in 
				\left]0,1\right].
			\end{equation}
			
			\item\label{t:genDRtau3} Let 
			$\theta\in\RPP$ be such that $(1-\theta)\in[-2\beta_{\cA},2\beta_{\cB}]$.
			Then 
			$R_{\theta \mathcal{B}}^{1/\theta}\circ R_{\mathcal{A}}^{\theta}$
			is $\alpha(\theta)$-averaged nonexpansive, where 
			$\alpha(\theta)\in\left]0,1\right]$ is given by
			\begin{equation}
				\label{e:avcompgentau}
				\alpha(\theta)=
				\begin{cases}
					1,&\text{if}\:\:\beta_{\cA}=\beta_{\cB}=0;\\
					\dfrac{2(1+\theta)(\beta_{\cA}+\beta_{\cB})}
					{4(\beta_{\cB}+\beta_{\cA}\beta_{\cB}+\theta\beta_{\cA})
						-(\theta-1)^2},&\text{if}\:\:\beta_{\cA}+\beta_{\cB}>0. 
				\end{cases}
			\end{equation}
			
		\end{enumerate}
	\end{proposition}
	\begin{proof}
		\ref{t:genDRtau1}: Let $\theta>0$ and let $x\in\HH$. We have
		\begin{align}
			0\in (\mathcal{A}+\mathcal{B})x\ \ 
			&\Leftrightarrow\ \ (\exists y\in\HH)\quad y- 
			x\in\mathcal{A}x\quad\text{and}\quad 
			x-y\in\mathcal{B}x\nonumber\\
			&\Leftrightarrow\ \ (\exists y\in\HH)\quad 
			x=J_{\mathcal{A}}y\quad\text{and}\quad 
			(1+\theta)x-\theta y\in(\Id+\theta\mathcal{B})x\nonumber\\
			&\Leftrightarrow\ \ (\exists y\in\HH)\quad 
			x=J_{\mathcal{A}}y\quad\text{and}\quad 
			x=J_{\theta\mathcal{B}}((1+\theta)x-\theta y)\nonumber\\
			&\Leftrightarrow\ \ (\exists y\in\HH)\quad 
			x=J_{\mathcal{A}}y\quad\text{and}\quad 
			x=J_{\theta\mathcal{B}}(R_{\mathcal{A}}^\theta y)\nonumber\\
			&\Leftrightarrow \ \ (\exists y\in\HH)\quad 
			x=J_{\mathcal{A}}y\quad\text{and}\quad 
			y=\left(1+\frac{1}{\theta}\right)x-\frac{1}{\theta}R_{\mathcal{A}}^\theta y=
			R_\mathcal{\theta 
				B}^{1/\theta}R_\mathcal{A}^{\theta}y\nonumber\\
			&\Leftrightarrow\ \ (\exists y\in\Fix R_\mathcal{\theta 
				B}^{1/\theta}R_\mathcal{A}^{\theta})\quad x=J_{\mathcal{A}}y.
		\end{align}
		
		\ref{t:genDRtau2}:
		Let $x$ and $y$ be in $\HH$. Since $J_{\mathcal{A}}$ is 
		$(2+2\beta_{\cA})^{-1}$-averaged nonexpansive 
		\eqref{e:equivcoco}, it 
		follows from 
		\cite[Corollary~2.15]{bauschke2011:Convex_analysis} that
		\begin{align}
			\|R_\mathcal{A}^{\theta}x-R_\mathcal{A}^{\theta}y\|^2
			&=(1+\theta)(\|J_\mathcal{A}x-J_\mathcal{A}y\|^2+\theta\|(\Id-J_{\mathcal{A}})x-(\Id-J_\mathcal{A})y\|^2)
			-\theta\|x-y\|^2\nonumber\\
			&\le 
			\|x-y\|^2-(1+\theta)(2\beta_{\cA}+1
			-\theta)\|(\Id-J_{\mathcal{A}})x-(\Id-J_\mathcal{A})y\|^2,
		\end{align}
		and \eqref{e:ineqrelref} follows.
		
		\ref{t:genDRtau2+}:
		Noting that 
		$\Id-R_{\mathcal{A}}^{\theta}=(1+\theta)(\Id-J_{\mathcal{A}})$, the 
		result follows from \ref{t:genDRtau2}.

		\ref{t:genDRtau3}: Since $\mathcal{B}$ is $\beta_{\cB}$-cocoercive and $\theta>0$,
		$\theta\mathcal{B}$ is $(\beta_{\cB}/\theta)$-cocoercive. Moreover, 
		\begin{equation}
			\label{e:equivcond}
			1-\theta\in[-2\beta_{\cA},2\beta_{\cB}]
			\quad\Leftrightarrow\quad
			\theta\le 1+2\beta_{\cA}
			\quad\text{and}\quad
			\frac{1}{\theta}\le 1+\frac{2\beta_{\cB}}{\theta},
		\end{equation}
		which are exactly the conditions of \ref{t:genDRtau2+} applied to
		$\mathcal{A}$ with parameter $\theta$ and to $\theta\mathcal{B}$ with
		parameter $1/\theta$. Hence $R_{\mathcal{A}}^{\theta}$ and
		$R_{\theta\mathcal{B}}^{1/\theta}$ are $\alpha_{\mathcal{A}}(\theta)$- and
		$\alpha_{\theta\mathcal{B}}(1/\theta)$-averaged nonexpansive, respectively,
		where
		\begin{equation}
			\label{e:twoaver}
			\alpha_{\mathcal{A}}(\theta)=\frac{1+\theta}{2(1+\beta_{\cA})}
			\qquad\text{and}\qquad
			\alpha_{\theta\mathcal{B}}(1/\theta)
			=\frac{1+\frac{1}{\theta}}{2+\frac{2\beta_{\cB}}{\theta}}
			=\frac{1+\theta}{2(\theta+\beta_{\cB})}.
		\end{equation}
		It follows from \eqref{e:equivcond} that
		$\alpha_{\mathcal{A}}(\theta)\le 1$ and
		$\alpha_{\theta\mathcal{B}}(1/\theta)\le 1$, with equality in the first case if
		and only if $\theta=1+2\beta_{\cA}$ and in the second case if and only if
		$\theta=1-2\beta_{\cB}$. 
		If $\beta_{\cA}=\beta_{\cB}=0$, then $(1-\theta)\in[-2\beta_{\cA},2\beta_{\cB}]$ forces $\theta=1$
		and both $R_{\mathcal{A}}^{1}$ and $R_{\mathcal{B}}^{1}$ are nonexpansive, and
		hence so is their composition, which is the case $\alpha(\theta)=1$. Suppose
		from now on that $\beta_{\cA}+\beta_{\cB}>0$.
		Hence,
		\cite[Proposition~2.4]{Yamada} asserts that
		$R_{\theta \mathcal{B}}^{1/\theta}R_{\mathcal{A}}^{\theta}$ is
		$\alpha(\theta)$-averaged nonexpansive, where
		\begin{equation}
			\label{e:yamada}
			\alpha(\theta)=\frac{\alpha_{\mathcal{A}}(\theta)+\alpha_{\theta\mathcal{B}}(1/\theta)
				-2\alpha_{\mathcal{A}}(\theta)\alpha_{\theta\mathcal{B}}(1/\theta)}
			{1-\alpha_{\mathcal{A}}(\theta)\alpha_{\theta\mathcal{B}}(1/\theta)}\in\left]0,1\right],
		\end{equation}
		noting that the proof of \cite[Proposition~2.4]{Yamada} is still valid if either $\alpha_{\mathcal{A}}(\theta)=1$ or $\alpha_{\theta\mathcal{B}}(1/\theta)=1$.
		Finally, \eqref{e:avcompgentau} follows from \eqref{e:yamada} and \eqref{e:twoaver} after
		simple computations.
	\end{proof}
	
	\begin{remark}\label{rem:DP}
		Since $R_{\mathcal{A}}^{1}$ is the classical 
		reflection, note that we recover \cite[Proposition~5.3]{Giselsson2017}
		from Proposition~\ref{t:genDRtau}.\ref{t:genDRtau2+}. 
	\end{remark}
	
	\begin{proposition}\label{t:RDRweakc}
		In the context of Problem~\ref{prob:main}, let $\tau\in\RPP$, let
		$\theta\in\RPP$, let $(\lambda_n)_{n\in\NN}$ be a sequence in $\RP$, and let
		$z_0\in\HH$. Consider the recurrence
		\begin{align}
			(\forall n\in\NN)\hspace{0.2cm}
			&\left\lfloor 
			\begin{array}{ll}
				x_{n}=J_{\tau \mathcal{A}}z_n\\
				y_n=J_{\theta\tau\mathcal{B}}\big((1+\theta)x_{n}
				-\theta z_n\big)\\
				z_{n+1}=z_n+\lambda_n\left(\dfrac{1+\theta}{\theta}\right)
				(y_n-x_{n}).
				\label{e:RDR}
			\end{array}
			\right.
		\end{align}
		Then, the following hold:
		\begin{enumerate}
			\item\label{p:sharpcoco0} 
			Suppose that $\lambda_n\in\left[0,1/\alpha\right]$ is such that 
			$\sum_{n\in\NN}\lambda_n(1-\alpha\lambda_n)=+\infty$, that $(1-\theta)\tau\in[-2\beta_{\cA},2\beta_{\cB}]$, where
			$\alpha\in\left]0,1\right]$ is defined by
			\begin{equation}
				\label{e:defalpha}
				\alpha=
				\begin{cases}
					1,&\text{if}\:\:\beta_{\cA}=\beta_{\cB}=0;\\
					\dfrac{2\tau(1+\theta)(\beta_{\cA}+\beta_{\cB})}
					{4(\tau\beta_{\cB}+\beta_{\cA}\beta_{\cB}+\tau\theta\beta_{\cA})
						-\tau^2(\theta-1)^2},&\text{if}\:\:\beta_{\cA}+\beta_{\cB}>0,    
				\end{cases}
			\end{equation}
			Then, there exists $\widehat{z}\in\Fix R_{\theta\tau 
				\mathcal{B}}^{1/\theta}R_{\tau\mathcal{A}}^{\theta}$ such that 
			$z_n\rightharpoonup\widehat{z}$ and $x_n\rightharpoonup 
			J_{\tau\mathcal{A}}\widehat{z}\in Z$.
			\item\label{p:sharpcoco1} Suppose that $\lambda_n\equiv1$. Then, $(z_n)_{n\in\NN}$ converges weakly for every real 
			Hilbert space $\HH$, every $\mathcal{A}$ and $\mathcal{B}$ satisfying the assumptions of
			Problem~\ref{prob:main}, and every $z_0\in\HH$ if and only if 
			$(1-\theta)\tau\in\left]-2\beta_{\cA},2\beta_{\cB}\right[$.
		\end{enumerate}
	\end{proposition}
	\begin{proof}
		\ref{p:sharpcoco0}:
		Note that $\tau\mathcal{A}$ and $\tau\mathcal{B}$
		are $(\beta_{\cA}/\tau)$- and $(\beta_{\cB}/\tau)$-cocoercive, 
		respectively, and that $Z=(\tau\mathcal{A}+\tau\mathcal{B})^{-1}(0)$. 
		Hence, it follows from \eqref{e:relaxedrefl} and simple computations 
		that
		\begin{equation}
			R_{\theta \tau\mathcal{B}}^{1/\theta}\circ 
			R_{\tau\mathcal{A}}^{\theta}
			=\Id+\left(\dfrac{1+\theta}{\theta}\right)
			\Big(J_{\theta\tau\mathcal{B}}\big((1+\theta)
			J_{\tau\mathcal{A}}-\theta\Id\big)-J_{\tau\mathcal{A}}\Big)
		\end{equation}
		and we deduce from \eqref{e:RDR} that 
		\begin{equation}
			(\forall n\in\NN)\quad z_{n+1}=(1-\lambda_n)z_n+\lambda_nTz_n,
		\end{equation}
		where $T=R_{\theta 
			\tau\mathcal{B}}^{1/\theta} 
		R_{\tau\mathcal{A}}^{\theta}$. In addition, 
		Proposition~\ref{t:genDRtau}.\ref{t:genDRtau3} implies that 
		$T$ is $\alpha$-averaged nonexpansive and
		\cite[Proposition~5.16]{bauschke2011:Convex_analysis} asserts that 
		there exists $\widehat{z}\in\Fix R_{\theta 
			\tau\mathcal{B}}^{1/\theta} 
		R_{\tau\mathcal{A}}^{\theta}$ such that 
		\begin{equation}
			\label{e:auxs}
			z_n\rightharpoonup 
			\widehat{z}\quad \text{and}\quad z_n-Tz_{n}\to 0.
		\end{equation}
		In view of Proposition~\ref{t:genDRtau}.\ref{t:genDRtau1}, we 
		obtain
		$J_{\tau\mathcal{A}}\widehat{z}\in Z$. Finally, by denoting 
		$\widehat{x}=J_{\tau\mathcal{A}}\widehat{z}$, the nonexpansivity 
		of $J_{\tau\mathcal{A}}$ yields
		$\|x_{n}-\widehat{x}\|\le \|z_n-\widehat{z}\|$ and, since 
		$(z_n)_{n\in\NN}$ is bounded, so is $(x_n)_{n\in\NN}$.
		Let $x$ be any weak accumulation point of $(x_n)_{n\in\NN}$, say 
		$x_{n_k}\rightharpoonup x$. Then 
		$z_{n_k}\rightharpoonup\widehat{z}$. Denote
		\begin{equation}
			(\forall k\in\NN)\quad 
			u_{n_k}=(z_{n_k}-x_{n_k})/\tau\quad\text{and}\quad 
			v_{n_k}=((1+\theta)x_{n_k}-\theta z_{n_k}-y_{n_k})/(\tau\theta).
		\end{equation}
		Hence, $(x_{n_k},u_{n_k})\in\gra\mathcal{A}$,
		$(y_{n_k},v_{n_k})\in\gra\mathcal{B}$, $u_{n_k}\rightharpoonup 
		(\widehat{z}-x)/\tau$, 
		$x_{n_k}-y_{n_k}=\theta(z_{n_k}-Tz_{n_k})/(1+\theta)\to 0$,
		and 
		$u_{n_k}+v_{n_k}=(x_{n_k}-y_{n_k})/\tau\theta\to 0$.
		Hence, it follows from 
		\cite[Corollary~26.6]{bauschke2011:Convex_analysis} that 
		$(\widehat{z}-x)/\tau\in\mathcal{A}x$ and, therefore, 
		$x=\widehat{x}=J_{\tau\mathcal{A}}\widehat{z}$.
		Therefore, there exists a unique weak accumulation point 
		of $(x_n)_{n\in\NN}$ and the result follows from \cite[Lemma~2.46]{bauschke2011:Convex_analysis}.
		First note that a direct computation yields
		\begin{equation}
			\label{e:sharpident}
			4(\tau\beta_{\cB}+\beta_{\cA}\beta_{\cB}+\tau\theta\beta_{\cA})
			-\tau^2(\theta-1)^2
			=\big(2\beta_{\cA}+(1-\theta)\tau\big)\big(2\beta_{\cB}-(1-\theta)\tau\big)
			+2\tau(1+\theta)(\beta_{\cA}+\beta_{\cB}).
		\end{equation}
		
		\ref{p:sharpcoco1}: Suppose that 
		$(1-\theta)\tau\in\left]-2\beta_{\cA},2\beta_{\cB}\right[$. Then 
		$\beta_{\cA}+\beta_{\cB}>0$ and \eqref{e:sharpident} implies that the 
		denominator of $\alpha$ defined in \eqref{e:defalpha} is positive and that
		\begin{equation*}
			1-\alpha=\frac{\big(2\beta_{\cA}+(1-\theta)\tau\big)
				\big(2\beta_{\cB}-(1-\theta)\tau\big)}
			{4(\tau\beta_{\cB}+\beta_{\cA}\beta_{\cB}+\tau\theta\beta_{\cA})
				-\tau^2(\theta-1)^2}>0.
		\end{equation*}
		Hence, $\lambda_n\equiv1$ satisfies $\lambda_n\in[0,1/\alpha]$ and 
		$\sum_{n\in\NN}\lambda_n(1-\alpha\lambda_n)=+\infty$, and the weak 
		convergence of $(z_n)_{n\in\NN}$ follows from 
		\ref{p:sharpcoco0}. Conversely, let $\HH=\RR$ and $z_0=1$, and 
		recall that
		\eqref{e:RDR} with $\lambda_n\equiv1$ reads 
		$z_{n+1}=R_{\theta\tau\mathcal{B}}^{1/\theta}(R_{\tau\mathcal{A}}^{\theta}z_n)$. 
		If $(1-\theta)\tau\le-2\beta_{\cA}$, set $\mathcal{B}=0$, and set 
		$\mathcal{A}=\beta_{\cA}^{-1}\Id$ if $\beta_{\cA}>0$ or 
		$\mathcal{A}=N_{\{0\}}$ if $\beta_{\cA}=0$. Then $\mathcal{A}$ is maximally 
		monotone and $\beta_{\cA}$-cocoercive, $Z=\{0\}$, 
		$J_{\tau\mathcal{A}}=\frac{\beta_{\cA}}{\beta_{\cA}+\tau}\Id$, and 
		$R_{\theta\tau\mathcal{B}}^{1/\theta}=\Id$, which yield, for every $n\in\NN$,
		$z_{n+1}=\frac{\beta_{\cA}-\theta\tau}{\beta_{\cA}+\tau}\,z_n$,
		where $(\beta_{\cA}-\theta\tau)/(\beta_{\cA}+\tau)\le-1$ since 
		$(1-\theta)\tau\le-2\beta_{\cA}$. Hence, $(z_n)_{n\in\NN}$ alternates its 
		sign and $|z_n|\ge1$ for every $n\in\NN$, so it does not converge. 
		Similarly, if $(1-\theta)\tau\ge2\beta_{\cB}$, set $\mathcal{A}=0$, and set 
		$\mathcal{B}=\beta_{\cB}^{-1}\Id$ if $\beta_{\cB}>0$ or 
		$\mathcal{B}=N_{\{0\}}$ if $\beta_{\cB}=0$. Then 
		$J_{\theta\tau\mathcal{B}}=\frac{\beta_{\cB}}{\beta_{\cB}+\theta\tau}\Id$ 
		and $R_{\tau\mathcal{A}}^{\theta}=\Id$, which yield 
		$z_{n+1}=\frac{\beta_{\cB}-\tau}{\beta_{\cB}+\theta\tau}\,z_n$, where 
		$(\beta_{\cB}-\tau)/(\beta_{\cB}+\theta\tau)\le-1$ since 
		$(1-\theta)\tau\ge2\beta_{\cB}$, and the same conclusion holds.
	\end{proof}
	
	\begin{remark}
		\begin{enumerate}
			\item In the classical case $\theta=1$, the weak convergence of the shadow
			sequence $(x_n)_{n\in\NN}$ goes back to \cite{Svai11}, and it was obtained
			for the algorithm considered here in \cite{BartzCampoyPhan2020} by means of
			demiclosedness principles for generalized nonexpansive mappings.
			
			\item 
			The operator $R_{\theta\tau\mathcal{B}}^{1/\theta}\circ
			R_{\tau\mathcal{A}}^{\theta}$ coincides with the adaptive Douglas--Rachford
			operator of \cite{DaoPhan2019} for the parameters $\lambda=1+\theta$,
			$\mu=1+1/\theta$, $\gamma=\tau$, and $\delta=\theta\tau$, which satisfy the
			coupling $(\lambda-1)(\mu-1)=1$ and $\delta=(\lambda-1)\gamma$ required
			there, and Proposition~\ref{t:genDRtau}.\ref{t:genDRtau1} is the
			corresponding specialization of \cite[Lemma~4.1(iii)]{DaoPhan2019}. The two
			convergence analyses are, however, driven by different regularity
			assumptions. In \cite{DaoPhan2019}, $\mathcal{A}$ and $\mathcal{B}$ are
			assumed maximally $a$- and $b$-monotone with $a+b\ge 0$, and the admissible
			parameters are those satisfying $\mu\in[2-2\gamma b,\,2+2\gamma a]$. Note that, in the
			setting of Problem~\ref{prob:main}, where $\mathcal{A}$ and $\mathcal{B}$
			are maximally monotone, this condition forces $a=b=0$ and hence $\theta=1$,
			that is, the classical Douglas--Rachford splitting. The range
			$(1-\theta)\in[-2\beta_{\cA},2\beta_{\cB}]$ obtained in
			Proposition~\ref{t:genDRtau}.\ref{t:genDRtau3} is instead a consequence of
			cocoercivity. Conversely, the analysis in \cite{DaoPhan2019} covers weakly
			monotone operators, which are excluded here since cocoercivity implies
			monotonicity. We also note that, in the absence of cocoercivity, the
			composition $R_{\theta\tau\mathcal{B}}^{1/\theta}\circ
			R_{\tau\mathcal{A}}^{\theta}$ fails to be averaged nonexpansive for
			$\theta\neq1$, even for subdifferential operators, and its convergence
			analysis in \cite{NilssonAkermanGiselsson2025} proceeds by a Lyapunov
			argument. Proposition~\ref{t:genDRtau}.\ref{t:genDRtau3} shows that
			cocoercivity restores averagedness, with the explicit constant
			$\alpha(\theta)$.
			
			\item 
			Since an operator is $\beta$-cocoercive if and only if it is
			$\beta$-comonotone with $\beta\ge0$, Problem~\ref{prob:main} falls within
			the scope of \cite[Theorem~3.3]{AlcantaraDaoTakeda2025}, which establishes
			the weak convergence of \eqref{e:RDR} with a constant relaxation
			$\lambda_n\equiv\kappa\in\left]0,1/\alpha\right[$ by combining the
			Attouch--Th\'era duality \cite{AttouchThera1996} with the conical averagedness framework of
			\cite{BartzDaoPhan2022,BauschkeMoursiWang2021}; the relaxation bound therein coincides with
			$1/\alpha$. When $\beta_{\cA}+\beta_{\cB}>0$, the admissible parameters in
			\cite{AlcantaraDaoTakeda2025} are determined by the inequality
			$\tau^2(1+\theta)^2<4(\tau+\beta_{\cA})(\theta\tau+\beta_{\cB})$, which strictly
			contains the interval $(1-\theta)\tau\in[-2\beta_{\cA},2\beta_{\cB}]$ of
			Proposition~\ref{t:RDRweakc}; outside this interval the composition
			$R_{\theta\tau\mathcal{B}}^{1/\theta}\circ R_{\tau\mathcal{A}}^{\theta}$
			need not be averaged nonexpansive. 
			Proposition~\ref{t:RDRweakc} provides, in
			turn, varying relaxations $(\lambda_n)_{n\in\NN}$ and a direct primal proof
			based on the averagedness properties established in
			Proposition~\ref{t:genDRtau}, which are also the starting point of the
			linear convergence analysis of Section~\ref{sec:strongmon}.
			We also note that, by
			\cite[Remark~3.5]{AlcantaraDaoTakeda2025}, the weak convergence of the
			shadow sequence cannot be improved to strong convergence under cocoercivity
			alone, which motivates assumption \eqref{e:strongmon} below.
		\end{enumerate}
		
	\end{remark}
	
	We now study in detail the cases when only one operator is 
	cocoercive. The following result is a consequence of 
	Proposition~\ref{t:RDRweakc}.\ref{p:sharpcoco0} by setting $\beta_{\cA}=0$ or $\beta_{\cB}=0$.
	\begin{corollary}
		In the context of Problem~\ref{prob:main}, let $\tau\in\RPP$, let 
		$\theta\in\RPP$,
		and suppose that one of the following holds:
		\begin{enumerate}
			\item $\beta_{\cA}=0$, $\beta_{\cB}>0$,
			$\theta\in 
			\left[\max\{0,1-2\beta_{\cB}/\tau\},1\right]$, 
			$(\lambda_n)_{n\in\NN}$ is a sequence in $\left[0,1/\overline{\alpha}_{\mathcal{B}}(\tau,\theta)\right]$ 
			such 
			that 
			$\sum_{n\in\NN}\lambda_n(1-\overline{\alpha}_{\mathcal{B}}(\tau,\theta)\lambda_n)=+\infty$, where 
			\begin{equation}
				\overline{\alpha}_{\mathcal{B}}(\tau,\theta)=\dfrac{2\beta_{\cB}(1+\theta)}
				{4\beta_{\cB}-\tau(\theta-1)^2}
				\in\left]0,1\right].
			\end{equation} 
			
			\item $\beta_{\cB}=0<\beta_{\cA}$,
			$\theta\in 
			\left[1,1+2\beta_{\cA}/\tau\right]$, $(\lambda_n)_{n\in\NN}$ is a sequence in $\left[0,1/\overline{\alpha}_{\mathcal{A}}(\tau,\theta)\right]$ 
			such 
			that 
			$\sum_{n\in\NN}\lambda_n(1-\overline{\alpha}_{\mathcal{A}}(\tau,\theta)\lambda_n)=+\infty$, where 
			\begin{equation}
				\overline{\alpha}_{\mathcal{A}}(\tau,\theta)=\dfrac{2\beta_{\cA}(1+\theta)}
				{4\theta\beta_{\cA}-\tau(\theta-1)^2}
				\in\left]0,1\right].
			\end{equation} 
		\end{enumerate}
		Then, the sequence $(z_n)_{n\in\NN}$ generated by \eqref{e:RDR}
		satisfies $z_n\rightharpoonup\widehat{z}$ and 
		$x_n\rightharpoonup 
		J_{\tau\mathcal{A}}\widehat{z}\in Z$ for some $\widehat{z}\in\Fix R_{\theta\tau 
			\mathcal{B}}^{1/\theta}R_{\tau\mathcal{A}}^{\theta}$.
	\end{corollary}
	
	\begin{remark}
		Note that, for small values of functions $\overline{\alpha}_{\mathcal{A}}$ and 
		$\overline{\alpha}_{\mathcal{B}}$ we obtain larger over-relaxation steps 
		$(\lambda_n)_{n\in\NN}$.  After simple computations, 
		we deduce the following:
		\begin{enumerate}
			\item In the case when $\beta_{\cA}=0<\beta_{\cB}$, we have that, for every 
			$\tau\in\left]0,4\beta_{\cB}/3\right]$,  
			$\overline{\alpha}_{\mathcal{B}}(\tau,\cdot)$ is non-decreasing. Therefore, 
			$\inf \overline{\alpha}_{\mathcal{B}}(\tau,\left]0,1\right])$ is obtained when 
			$\theta\to 0$. However, $\theta=0$ cannot be used in 
			\eqref{e:RDR}.
			Now if $\tau>4\beta_{\cB}/3$, $\overline{\alpha}_{\mathcal{B}}(\tau,\cdot)$ is 
			minimized at $\theta_{\mathcal{B}}(\tau)=-1+2\sqrt{1-\beta_{\cB}/\tau}$, 
			leading to 
			$$(\forall \tau>4\beta_{\cB}/3)\quad 
			\overline{\alpha}_{\mathcal{B}}(\tau,\theta_{\mathcal{B}}(\tau))
			=\frac{1+\sqrt{1-\beta_{\cB}/\tau}}{2}.$$
			Note that, as $\tau\to 4\beta_{\cB}/3$, 
			$\theta_{\mathcal{B}}(\tau)\to 0$ and 
			$\overline{\alpha}_{\mathcal{B}}(\tau,\theta_{\mathcal{B}}(\tau))\to 3/4$.
			
			\item In the case when $\beta_{\cB}=0<\beta_{\cA}$, $\overline{\alpha}_{\mathcal{A}}(\tau,\cdot)$
			is minimized at 
			$\theta_{\mathcal{A}}(\tau)=-1+2\sqrt{1+\beta_{\cA}/\tau}$, leading to 
			$$(\forall \tau>0)\quad 
			\overline{\alpha}_{\mathcal{A}}(\tau,\theta_{\mathcal{A}}(\tau))
			=\frac{1+\frac{1}{\sqrt{1+\beta_{\cA}/\tau}}}{2}.$$
			Note that, as $\tau\to 0$, $\theta_{\mathcal{A}}(\tau)\to +\infty$ and 
			$\overline{\alpha}_{\mathcal{A}}(\tau,\theta_{\mathcal{A}}(\tau))\to 1/2$.
		\end{enumerate}
		Therefore, in the case when only one operator is cocoercive, say 
		$\mathcal{A}$ is $\beta_{\cA}$-cocoercive for some $\beta_{\cA}>0$, the 
		algorithm in \eqref{e:RDR} admits larger relaxation steps if the 
		step-size $\tau$ is chosen small, 
		$\theta=-1+2\sqrt{1+\beta_{\cA}/\tau}$, and 
		$\lambda_n\approx 1/\alpha$, where 
		$\alpha=(1+\sqrt{1+\beta_{\cA}/\tau})/(2\sqrt{1+\beta_{\cA}/\tau})$.
		Note that the cocoercive operator is firstly implemented at each 
		iteration in 
		\eqref{e:RDR}.
	\end{remark}
	
	\begin{remark} Note that, since $\beta_{\cA}\ge0$ and $\beta_{\cB}\ge 0$ in Problem~\ref{prob:main}, the convergence of the sequences generated by \eqref{e:RDR}
		is guaranteed by 
		choosing $\theta=1\in 
		\left[\max\{0,1-2\beta_{\cB}/\tau\},1+2\beta_{\cA}/\tau\right]\smallsetminus\{0\}$.
		In this case, 
		$\alpha=\tau(\beta_{\cA}+\beta_{\cB})/(\tau(\beta_{\cA}+\beta_{\cB})+\beta_{\cA}\beta_{\cB})\in\left]0,1\right]$ and
		\eqref{e:RDR} reduces to the Douglas-Rachford splitting with a 
		relaxation step. Observe that, when $\mathcal{A}$ or $\mathcal{B}$ is merely maximally monotone,
		($\beta_{\cA}=0$ or $\beta_{\cB}=0$) we have $\alpha=1$ 
		and the only choice guaranteeing unconditional convergence is $\theta=1$. This result is coherent with \cite[Theorem~3.2]{NilssonAkermanGiselsson2025}, which builds
		on the uniqueness theorem of \cite{Ryu2020UniquenessDRS}, unconditional convergence of \eqref{e:RDR}
		over the class of maximally monotone pairs forces $\theta=1$ and
		$\lambda_n(1+\theta)/\theta\in\left]0,2\right[$ whenever $\dim\HH\ge2$.
		For subdifferentials of functions in $\Gamma_0(\HH)$,
		\cite[Theorem~3.3]{NilssonAkermanGiselsson2025} shows that every
		$\theta\in\RPP$ is admissible without cocoercivity, provided
		$\lambda_n<2\min\{\theta,1\}/(1+\theta)$, and that this region is sharp.
	\end{remark}

	\subsection{Strongly monotone inclusions}
	\label{sec:strongmon}
	In the context of Problem~\ref{prob:main}, suppose 
	\begin{equation}
		\label{e:strongmon}\tag{$H$}
		\mathcal{A}\quad \text{is $\rho$-strongly monotone, for some 
			$\rho> 0$.}
	\end{equation}
	In the case when $\beta_{\cA}>0$, we have from \eqref{e:defcoco} and Cauchy-Schwarz that $\rho\in\left]0,\beta_{\cA}^{-1}\right]$.
	
	Since an operator is strongly monotone cocoercive if and only if its resolvent is averaged nonexpansive and cocoercive in view of \eqref{e:equivcoco} and \eqref{e:equivstrong}, we first provide a tight rate for conical relaxations of averaged nonexpansive and cocoercive operators.
	The following result is a generalization of
	\cite[Theorem~7.2]{Giselsson2017} and is interesting in its own right. 
	We note that the composition results of
	\cite{Giselsson2017} have been extended in a different direction, to
	compositions of conically nonexpansive operators, in
	\cite{GiselssonMoursi2021}.
	
	\begin{proposition}
		\label{p:reflectcontr}
		Let $\alpha\in\left]0,\frac{1}{2}\right]$, let $\beta\in\left]1,+\infty\right[$, and let 
		$T\colon\HH\to\HH$ be an $\alpha$-averaged nonexpansive and
		$\beta$-cocoercive operator. Then, we have
		$\beta(1-2\alpha)\le 1$ and the following assertions hold:
		\begin{enumerate}
			\item\label{p:reflectcontr1} For every $\lambda\in\left[0,+\infty\right[$, 
			the operator $S_{\lambda}=(1-\lambda)\Id+\lambda T$ is 
			$\rate_T(\lambda)$-Lipschitz continuous, where 
			{\small\begin{equation}
					\label{e:defrho}
					\vspace{-.2cm}
					\rate_T(\lambda)=
					\begin{cases}
						1-\lambda\left(1-\beta^{-1}\right), &\text{if}\:\:0\le\lambda<\dfrac{2\beta}{2\beta-1},\\[3mm]
						\sqrt{1-\lambda(1-\alpha\lambda)\bigg(
							\dfrac{2(\beta-1)}{2\beta(1-\alpha)-1}\bigg)}, &\text{if}\:\:\dfrac{2\beta}{2\beta-1}\le\lambda\le\dfrac{1}{\alpha},\\[3mm]
						2\alpha\lambda-1, &\text{if}\:\:\lambda>\dfrac{1}{\alpha}.
					\end{cases}
			\end{equation}}
			Moreover, by setting $\widehat{\lambda}=\max\{\frac{1}{2\alpha},\frac{2\beta}{2\beta-1}\}$, we have
			{\small\begin{equation*}
					\rate_T(\widehat{\lambda})=\min_{\lambda \in \left[0,+\infty\right[}\rate_T(\lambda)=
					\begin{cases}
						\sqrt{1-\dfrac{\beta-1}{2\alpha(2\beta(1-\alpha)-1)}}, &\text{if}\:\:\:\alpha\in \left]0, \frac{2\beta-1}{4\beta}\right[;\\[3mm]
						\dfrac{1}{2\beta-1},\quad &\text{if}\:\:\:\alpha\in \left[\frac{2\beta-1}{4\beta},\frac{1}{2}\right].\\
					\end{cases}
			\end{equation*}}
			
			\item \label{p:reflectcontr2} Given $\lambda\in\left[0,+\infty\right[$,
			the rate $\rate_T(\lambda)$ in \eqref{e:defrho} is attained by taking $\HH=H\oplus H$, for a real 
			Hilbert space $H$, and
			\begin{equation}
				\label{e:tightT}
				T\colon \HH \to \HH \colon
				(x,y) \mapsto 
				\begin{pmatrix}
					\mu x+\sigma y\\
					\mu y-\sigma x
				\end{pmatrix},
			\end{equation}
			where{\small
				\begin{equation} 
					\label{e:tightmusigma}
					\mu=
					\begin{cases}
						\dfrac{1}{\beta},&\hspace{-0.2cm}\text{if}\:\:\lambda<\dfrac{2\beta}{2\beta-1};\\[3mm]
						\dfrac{\beta(1-2\alpha)}{2\beta(1-\alpha)-1},&\hspace{-0.2cm}\text{if}\:\:\dfrac{2\beta}{2\beta-1}\le\lambda\le\dfrac{1}{\alpha};\\[3mm]
						1-2\alpha,&\hspace{-0.2cm}\text{if}\:\:\lambda>\dfrac{1}{\alpha}
					\end{cases}
					\ 
					\sigma=
					\begin{cases}
						\sqrt{\dfrac{\mu}{\beta}-\mu^2},&\hspace{-0.2cm}\text{if}\:\:\dfrac{2\beta}{2\beta-1}\le\lambda\le\dfrac{1}{\alpha};\\[3mm]
						0,&\text{otherwise}.
					\end{cases}
			\end{equation}}
			Indeed, $T$ is $\alpha$-averaged nonexpansive and $\beta$-cocoercive, and 
			$\|S_{\lambda}x-S_{\lambda}y\|=\rate_T(\lambda)\|x-y\|$ for every $x$ and $y$ 
			in $\HH$.
		\end{enumerate}
	\end{proposition}

	\begin{proof} 
		Let $x$ and $y$ be in $\HH$.  By expanding 
		$\|(\Id-T)x-(\Id-T)y\|^2$ in \eqref{e:averaged}, the 
		$\alpha$-averagedness of $T$ is equivalent to
		\begin{equation}
			\label{e:averaux}
			\|Tx-Ty\|^{2}-2(1-\alpha)\scal{x-y}{Tx-Ty}+(1-2\alpha)\|x-y\|^{2}\le 0.
		\end{equation}
		The $\beta$-cocoercivity of $T$ reads
		\begin{equation}
			\label{e:coco}
			\scal{x-y}{Tx-Ty}\ge\beta\|Tx-Ty\|^2
		\end{equation}
		and implies that $T$ is $\beta^{-1}$-Lipschitz continuous by the 
		Cauchy--Schwarz inequality. Moreover, by \cite[Proposition~4.35]{bauschke2011:Convex_analysis}, $T+(\alpha-1)\Id$ is $\alpha$-Lipschitz. 
		Hence, since $\scal{x-y}{Tx-Ty}=\scal{x-y}{Tx-Ty+(\alpha-1)(x-y)}+(1-\alpha)\|x-y\|^2$,
		the 
		Cauchy--Schwarz inequality and the $\beta^{-1}$-Lipschitz continuity of 
		$T$ yield
		\begin{equation}
			\label{e:sandwich}
			(1-2\alpha)\|x-y\|^2\le\scal{x-y}{Tx-Ty} \le\beta^{-1}\|x-y\|^2,
		\end{equation}
		which gives $\beta(1-2\alpha)\le1$ by taking $x\neq y$. Now let $\lambda\in\left[0,+\infty\right[$. We have
		\begin{equation}
			\label{e:tight0}
			\|S_{\lambda}x-S_{\lambda}y\|^2=\lambda^2\|Tx-Ty\|^2-2\lambda(\lambda-1)
			\scal{x-y}{Tx-Ty}+(\lambda-1)^2\|x-y\|^2
		\end{equation}
		and, by bounding $\lambda^2\|Tx-Ty\|^2$ via \eqref{e:coco} and via 
		\eqref{e:averaux}, respectively, we obtain
		\begin{align}
			\label{e:aux12a}
			\|S_{\lambda}x-S_{\lambda}y\|^2
			&\le\lambda\left(\lambda\beta^{-1}-2(\lambda-1)\right)
			\scal{x-y}{Tx-Ty}+(\lambda-1)^2\|x-y\|^2,\\
			\label{e:aux12}
			\|S_{\lambda}x-S_{\lambda}y\|^2
			&\le 2\lambda(1-\lambda\alpha)\scal{x-y}{Tx-Ty}
			+\big((\lambda-1)^2-\lambda^2(1-2\alpha)\big)\|x-y\|^2.
		\end{align}
		The coefficient of $\scal{x-y}{Tx-Ty}$ in \eqref{e:aux12a} is nonnegative 
		if and only if $\lambda\le2\beta/(2\beta-1)$, and the one in 
		\eqref{e:aux12} is nonnegative if and only if $\lambda\le1/\alpha$. Note 
		that, since $\alpha\le \frac{1}{2}$, we have
		\begin{equation}
			\label{e:condpos1}
			2\beta(1-\alpha)\ge\beta>1, \quad
			\textnormal{
				which implies} \quad
			\frac{2\beta}{2\beta-1}<\frac{1}{\alpha}.
		\end{equation}
		We consider three cases to obtain \ref{p:reflectcontr1} and \ref{p:reflectcontr2}.
		
		\emph{Case 1: $0\le\lambda<2\beta/(2\beta-1)$.} Since the coefficient of the inner product term in \eqref{e:aux12a} is positive,
		the upper bound in \eqref{e:sandwich} implies that
		\begin{equation}
			\label{e:ineqcase1}
			\|S_{\lambda}x-S_{\lambda}y\|^2
			\leq\left(1-\lambda\left(1-\beta^{-1}\right)\right)^2\|x-y\|^2,
		\end{equation}
		where $1-\lambda(1-\beta^{-1})>0$, since 
		$\lambda(1-\beta^{-1})<2(\beta-1)/(2\beta-1)<1$. This gives the first case in \eqref{e:defrho}. The equality in \eqref{e:ineqcase1} is attained 
		in this case by $T=\beta^{-1}\Id$, since it satisfies  \eqref{e:coco} and the upper bound in \eqref{e:sandwich} with equality.
		
		\emph{Case 2: $2\beta/(2\beta-1)\le\lambda\le1/\alpha$.} Let 
		$\delta\in\left[1-2\alpha,\beta^{-1}\right]$. If $\scal{x-y}{Tx-Ty}\ge\delta\|x-y\|^2$, since the 
		coefficient in \eqref{e:aux12a} is nonpositive, we obtain
		\begin{equation*}
			\|S_{\lambda}x-S_{\lambda}y\|^2\le\big((\lambda-1)^2
			+\lambda(\lambda\beta^{-1}-2(\lambda-1))\delta\big)\|x-y\|^2.
		\end{equation*}
		Otherwise, since the coefficient in \eqref{e:aux12} is nonnegative, we 
		obtain
		\begin{equation}
			\label{e:ineqcase2}
			\|S_{\lambda}x-S_{\lambda}y\|^2\le\big(2\lambda(1-\lambda\alpha)\delta
			+(\lambda-1)^2-\lambda^2(1-2\alpha)\big)\|x-y\|^2.
		\end{equation}
		In order to find the smallest Lipschitz constant for $S_{\lambda}$, we 
		impose the equality of both constants and we obtain 
		$\delta^*=\beta(1-2\alpha)/(2\beta(1-\alpha)-1)\ge 0$. Noting that, in this case,
		\begin{equation*}
			(\lambda-1)^2+\lambda(\lambda\beta^{-1}-2(\lambda-1))
			\delta^*=1-\lambda(1-\alpha\lambda)\bigg(
			\dfrac{2(\beta-1)}{2\beta(1-\alpha)-1}\bigg)=\rate_T(\lambda)^2,
		\end{equation*}
		the second case in \eqref{e:defrho} holds. For the tightness, it follows from \eqref{e:tightmusigma} 
		that, in this case, 
		$\mu=\beta(1-2\alpha)/(2\beta(1-\alpha)-1)\ge0$ and
		\begin{equation*}
			\beta^{-1}-\mu=\frac{(\beta-1)\big(1-\beta(1-2\alpha)\big)}{\beta(2\beta(1-\alpha)-1)}\ge0,
		\end{equation*}
		in view of \eqref{e:condpos1} and since $\beta(1-2\alpha)\le 1$.
		Thus, $\mu/\beta-\mu^2=\mu(\beta^{-1}-\mu)\ge0$ and $\sigma$ is well defined. 
		Now let $J\colon\HH\to\HH\colon(x,y)\mapsto(y,-x)$, which is linear and 
		satisfies $\scal{z}{Jz}=0$ and $\|Jz\|=\|z\|$ for every $z\in\HH$. Noting that, in this case, 
		$T=\mu\Id+\sigma J$, $T$ is linear and
		\begin{equation}
			\label{e:tightTz}
			(\forall z\in\HH)\quad
			\scal{z}{Tz}=\mu\|z\|^2
			\quad\text{and}\quad
			\|Tz\|^2=(\mu^2+\sigma^2)\|z\|^2=\frac{\mu}{\beta}\|z\|^2.
		\end{equation}
		Therefore, the linearity of $T$ implies that \eqref{e:coco} and \eqref{e:averaux} hold with equality and, thus, \eqref{e:aux12a} and \eqref{e:aux12} too. Since $\mu=\delta^*$,
		the result follows from \eqref{e:tightTz}.
		
		\emph{Case 3: $\lambda>1/\alpha$.} Since the coefficient in 
		\eqref{e:aux12} is negative, the lower bound in \eqref{e:sandwich} yields
		\begin{equation}
			\label{e:ineqcase3}
			\|S_{\lambda}x-S_{\lambda}y\|^2
			\leq (2\alpha\lambda-1)^2\|x-y\|^2,
		\end{equation}
		which yields the third case in \eqref{e:defrho}. The equality in \eqref{e:ineqcase3} is obtained by taking $T=(1-2\alpha)\Id$, because it satisfies the lower bound in \eqref{e:sandwich} and \eqref{e:averaux} also with equality. Note that this case uses only the $\alpha$-averagedness of $T$ and, hence, it is valid without any cocoercivity assumption.

		Altogether, $S_{\lambda}$ is $\rate_T(\lambda)$-Lipschitz continuous for 
		every $\lambda\in\left[0,+\infty\right[$, where $\rate_T$ is defined in 
		\eqref{e:defrho} and the rate is attained by $T$ defined in \eqref{e:tightT}. 
		
		To compute the minimum of $\rate_T$, define 
		$\psi_{\alpha,\beta}\colon \lambda\mapsto 1-2\lambda(1-\lambda\alpha)
		(\beta-1)/(2\beta(1-\alpha)-1)$, so that $\rate_T^2=\psi_{\alpha,\beta}$ on 
		$[2\beta/(2\beta-1),1/\alpha]$. We have 
		$\psi_{\alpha,\beta}([2\beta/(2\beta-1),1/\alpha])\subset\left]-\infty,1\right]$ 
		and $\psi_{\alpha,\beta}(1/\alpha)=1$. In addition, $\beta(1-2\alpha)\leq 1$ 
		is equivalent to $2\alpha(2\beta(1-\alpha)-1)\ge \beta-1$ for 
		$\alpha \in \left]0,\frac{1}{2}\right[$. Hence, since 
		$\psi_{\alpha,\beta}$ is quadratic and convex, we deduce that
		\begin{equation*}
			\min_{\lambda \in \RR}\psi_{\alpha,\beta}(\lambda)= \psi_{\alpha,\beta}(1/(2\alpha))=1-\dfrac{\beta-1}{2\alpha(2\beta(1-\alpha)-1)}\ge 0
		\end{equation*}
		and that $\psi_{\alpha,\beta}$ is increasing in 
		$\left[1/(2\alpha),+\infty\right[$. Therefore, if 
		$\alpha\ge (2\beta-1)/(4\beta)$, we have 
		$1/(2\alpha)\le 2\beta/(2\beta-1)$ and, thus, 
		$\displaystyle \min_{\lambda \in [2\beta/(2\beta-1),1/\alpha]}\psi_{\alpha,\beta}(\lambda)= \psi_{\alpha,\beta}(2\beta/(2\beta-1))=1/(2\beta-1)^2$. 
		Finally, since $\beta>1$, $\rate_T$ is decreasing on 
		$\left[0,2\beta/(2\beta-1)\right[$ and increasing on 
		$\left]1/\alpha,+\infty\right[$. Moreover, $\rate_T$ is continuous, since 
		the expressions in \eqref{e:defrho} coincide at $\lambda=2\beta/(2\beta-1)$ 
		with value $1/(2\beta-1)$ and at $\lambda=1/\alpha$ with value $1$. 
		Therefore, the minimum of $\rate_T$ over $\left[0,+\infty\right[$ coincides 
		with its minimum over $[2\beta/(2\beta-1),1/\alpha]$, and the result 
		follows.
	\end{proof}
	The following result provides the contractive property of the 
	relaxed resolvent, which yields the convergence of 
	the algorithm with a tight linear rate.
	The worst case leading to the tightness result
	is inspired by the example provided in \cite[Section~7.1]{Giselsson2017}.
	\begin{theorem}\label{t:rrcontractive}
		In the context of Problem~\ref{prob:main}, assume 
		\eqref{e:strongmon}. Then, the following assertions hold.
		\begin{enumerate}
			\item\label{t:rrcontractive1}  For every 
			$\theta\in\left[-1,+\infty\right[$, $R_{\mathcal{A}}^{\theta}$  is 
			$\rate_{\mathcal{A}}(\theta)$-Lipschitz 
			continuous, where 
			\begin{equation}
				\label{e:defrhoa}
				\rate_{\mathcal{A}}(\theta)=
				\begin{cases}
					\dfrac{1-\rho\theta}{1+\rho}, &\text{if}\:\:-1\le\theta<\dfrac{1}{1+2\rho},\\[3mm]
					\sqrt{\dfrac{\beta_{\cA}+
							\rho\theta(\theta-2\beta_{\cA})}{\beta_{\cA}+\rho+2\rho\beta_{\cA}}},
					&\text{if}\:\:\dfrac{1}{1+2\rho}\le\theta\le 1+2\beta_{\cA},\\[3mm]
					\dfrac{\theta-\beta_{\cA}}{1+\beta_{\cA}}, &\text{if}\:\:\theta>1+2\beta_{\cA}.
				\end{cases}
			\end{equation}
			Moreover, $\rate_{\mathcal{A}}$ is minimized at 
			$\widehat{\theta}=\max\{(1+2\rho)^{-1},\beta_{\cA}\}$ and
			\begin{equation}
				\rate_{\mathcal{A}}(\widehat{\theta})=\min_{\theta\in\left[-1,+\infty\right[}\rate_{\mathcal{A}}(\theta)=
				\begin{cases}
					\rate_{\mathcal{A}}((1+2\rho)^{-1})=\dfrac{1}{1+2\rho},\quad 
					&\text{if}\:\:\beta_{\cA}\le \frac{1}{1+2\rho};\\[3mm]
					\rate_{\mathcal{A}}(\beta_{\cA})=\sqrt{\dfrac{\beta_{\cA}
							(1-\rho\beta_{\cA})}{\beta_{\cA}+\rho+2\rho\beta_{\cA}}},
					&\text{if}\:\:\beta_{\cA}> \frac{1}{1+2\rho}.
			\end{cases}	\end{equation} 
			\item\label{t:rrcontractive2}  Suppose that $\beta_{\cA}+\beta_{\cB}>0$, let $\tau>0$, and define
			\begin{equation}\label{e:I(tau)}
				I(\tau)=\left]0,1+\dfrac{2\beta_{\cA}}{\tau}\right[\cap\left[1-\dfrac{2\beta_{\cB}}{\tau},+\infty\right[.
			\end{equation}
			For every $\theta\in I(\tau)$,
			$R_{\theta \tau\mathcal{B}}^{1/\theta}\circ 
			R_{\tau\mathcal{A}}^{\theta}$ is
			$\rate(\theta,\tau)$-Lipschitz 
			continuous, 
			where \linebreak $\rate(\theta,\tau)\in\left[0,1\right[$ is given by
			\begin{equation}
				\label{e:defomega}
				\rate(\theta,\tau)=
				\begin{cases}
					\dfrac{1-\tau\rho\theta}{1+\tau\rho}, &\text{if}\:\:\theta<\dfrac{1}{1+2\tau\rho},\\[3mm]
					\sqrt{\dfrac{\beta_{\cA}+\tau\rho\theta(\tau
							\theta-2\beta_{\cA})}{\beta_{\cA}+\tau^2\rho+2\tau\beta_{\cA}\rho}}, &\text{if}\:\:\theta\ge\dfrac{1}{1+2\tau\rho}.
				\end{cases}
			\end{equation}
			Moreover, by setting 
			\begin{equation}\label{e:hatalpha}
				\widehat{\beta}=\max\left\{\beta_{\cA},\;\frac{\beta_{\cB}}{1+2\beta_{\cB}\rho}\right\},
			\end{equation}
			we have
			\begin{equation}\label{e:main}
				\rate(\widehat{\theta},\widehat{\tau})=
				\min_{\tau\in\left]0,+\infty\right[}\;
				\min_{\theta\in I(\tau)}\;\rate(\theta,\tau)
				=\frac{\sqrt{\widehat{\beta}\big(\frac{1}{\rho}-\widehat{\beta}\big)}}{\sqrt{\beta_{\cB}^2+2\widehat{\beta}\beta_{\cB}+\frac{\widehat{\beta}}{\rho}}+\widehat{\beta}+\beta_{\cB}},
			\end{equation}
			where 
			\begin{equation}\label{e:argmin}
				\widehat{\tau}=\beta_{\cB}+\sqrt{\beta_{\cB}^2+2\widehat{\beta}\beta_{\cB}+\frac{\widehat{\beta}}{\rho}}\in\RPP
				\qquad\text{and}\qquad
				\widehat{\theta}=1-\frac{2\beta_{\cB}}{\widehat{\tau}}\in\left]0,1\right].
			\end{equation}
			
			\item\label{t:rrcontractive3} 
			Let $\tau>0$ and $\theta\in\left[-1,+\infty\right[$, and suppose that 
			$\beta_{\cA}>0$ or $\theta<(1+2\tau\rho)^{-1}$. The constant
			$\rate_{\tau\mathcal{A}}(\theta)$ in \ref{t:rrcontractive1}, which coincides 
			with $\rate(\theta,\tau)$ in \eqref{e:defomega} if $\theta\in I(\tau)$, is 
			attained by taking $\HH=H\oplus H$, for a real Hilbert space $H$, 
			$\mathcal{B}=0$, and 
			\begin{equation}
				\mathcal{A}\colon \HH \to \HH \colon
				(x,y) \mapsto 
				\begin{pmatrix}
					\eta x-\xi y\\
					\eta y+\xi x
				\end{pmatrix},
			\end{equation}
			where
			\begin{equation}\label{eq:defmu}
				\eta=
				\begin{cases}
					\rho,&\hspace{-0.27cm}\text{if}\:\:\theta\le 1+\dfrac{2\beta_{\cA}}{\tau};\\[3mm]
					\beta_{\cA}^{-1},&\hspace{-0.27cm}\text{if}\:\:\theta> 1+\dfrac{2\beta_{\cA}}{\tau},
				\end{cases}
				\quad
				\xi=
				\begin{cases}
					\sqrt{\dfrac{\rho}{\beta_{\cA}}-\rho^2},\hspace{-0.27cm}&\text{if}\:\:\dfrac{1}{1+2\tau\rho}\le\theta\le 1+\dfrac{2\beta_{\cA}}{\tau};\\[3mm]
					0,&\text{otherwise}.
				\end{cases}
			\end{equation}
			
		\end{enumerate}
	\end{theorem}
	\begin{proof} \ref{t:rrcontractive1}: It follows from \eqref{e:equivcoco} and \eqref{e:equivstrong} that $J_{\mathcal{A}}$ is 
		$(1+\rho)$-cocoercive and $(2+2\beta_{\cA})^{-1}$-averaged 
		nonexpansive. Then, by setting $\alpha=(2+2\beta_{\cA})^{-1}\in\left]0,\frac{1}{2}\right]$ and 
		$\beta=1+\rho>1$, 
		Proposition~\ref{p:reflectcontr} applied to $\lambda = 1+\theta$ and $T = J_{\mathcal{A}}$ yields
		that $R_{\mathcal{A}}^{\theta}=S_{1+\theta}$ is $\rate_{J_{\mathcal{A}}}(1+\theta)$-Lipschitz continuous for every 
		$\theta\in\left[-1,+\infty\right[$. Since 
		$2\beta/(2\beta-1)=1+(1+2\rho)^{-1}$ and $1/\alpha=1+(1+2\beta_{\cA})$, 
		a direct computation shows that $\rate_{J_{\mathcal{A}}}(1+\theta)=\rate_{\mathcal{A}}(\theta)$, 
		where $\rate_{\mathcal{A}}$ is defined in \eqref{e:defrhoa}. Moreover, 
		$\widehat{\lambda}-1=\max\{\beta_{\cA},(1+2\rho)^{-1}\}=\widehat{\theta}$, 
		and the result follows. Note that the third case in \eqref{e:defrhoa} only uses the cocoercivity property of $\cA$ and, hence, it is valid without any strong monotonicity assumption.
		
		\ref{t:rrcontractive2}: Since $\theta\tau\mathcal{B}$ is $\beta_{\cB}/(\theta\tau)$-cocoercive for every $\theta>0$,
		if $\theta\ge 1-2\beta_{\cB}/\tau$ we deduce from 
		Proposition~\ref{t:genDRtau}.\ref{t:genDRtau2+} that 
		$R_{\theta\tau\mathcal{B}}^{1/\theta}$ is nonexpansive.  Moreover,
		$\tau\mathcal{A}$ is $\tau\rho$-strongly monotone and 
		$\beta_{\cA}/\tau$-cocoercive.
		Therefore, observing that 
		$\rate(\theta,\tau)=\rate_{\tau\mathcal{A}}(\theta)$ for every $\theta\in I(\tau)$, we deduce 
		from \ref{t:rrcontractive1} that
		\begin{equation}
			(\forall\theta\in I(\tau))
			\quad R_{\theta\tau\mathcal{B}}^{1/\theta}\circ 
			R_{\tau\mathcal{A}}^{\theta}\quad\text{is 
				$\rate(\theta,\tau)$-Lipschitz continuous.}
		\end{equation}
		Now, to optimize $\rate$, set 
		$I_0(\tau)=I(\tau)\cap\left[(1+2\tau\rho)^{-1},+\infty\right[$. If 
		$\theta\in I(\tau)\smallsetminus I_0(\tau)$, then 
		$1-2\beta_{\cB}/\tau\le\theta<(1+2\tau\rho)^{-1}<1+2\beta_{\cA}/\tau$, 
		which yields $(1+2\tau\rho)^{-1}\in I_0(\tau)$. Moreover, 
		$\rate(\cdot,\tau)$ is decreasing on $\left]0,(1+2\tau\rho)^{-1}\right]$, 
		since both expressions in \eqref{e:defomega} coincide at 
		$\theta=(1+2\tau\rho)^{-1}$ with value $(1+2\tau\rho)^{-1}$. Hence, 
		$\rate(\theta,\tau)>\rate((1+2\tau\rho)^{-1},\tau)$, and the minimum of 
		$\rate(\cdot,\tau)$ over $I(\tau)$ coincides with its minimum over 
		$I_0(\tau)$, on which $\rate(\cdot,\tau)$ is given by the second 
		expression in \eqref{e:defomega}. We consider the change of variables 
		$s=\theta\tau$, which yields
		{\small\begin{equation}
				\theta\in I_0(\tau)\quad\Leftrightarrow\quad s\in \tau I_0(\tau)=\left[\max\left\{\tau-2\beta_{\cB},\frac{\tau}{1+2\tau\rho}\right\},\tau+2\beta_{\cA}\right[,
		\end{equation}}
		amd, thus,
		\begin{equation}
			\label{e:defr}
			r(s,\tau):=\rate(s/\tau,\tau)
			=\sqrt{\dfrac{\beta_{\cA}+\rho s(s-2\beta_{\cA})}{\beta_{\cA}+\tau^2\rho+2\tau\beta_{\cA}\rho}}.
		\end{equation}
		Now, since it is clear that
		\begin{equation}
			\min_{\tau\in\left]0,+\infty\right[}\;
			\min_{\theta\in I(\tau)}\;\rate(\theta,\tau)=
			\min_{\tau\in\left]0,+\infty\right[}\;
			\min_{s\in \tau I_0(\tau)}\;r(s,\tau)=:\widehat{r},
		\end{equation}
		we solve the simpler separable last problem.
		Note that the denominator in \eqref{e:defr} is increasing in $\tau$ and
		{\small\begin{equation}
				(\forall s>0)\quad s\in \tau I_0(\tau)\quad\Leftrightarrow\quad 
				\tau\in \Sigma(s):=\left]s-2\beta_{\cA}, \min\left\{s+2\beta_{\cB},\left(\frac{s}{1-2s\rho}\right)_{\!\!+}\right\}\right].
		\end{equation}}
		Hence, by defining 
		\begin{equation}
			\label{e:defsstar}
			s^*=\beta_{\cB}\left(\sqrt{1+\frac{1}{\beta_{\cB}\rho}}-1\right)\in\left[0,\frac{1}{2\rho}\right[,
		\end{equation} 
		which is where functions $s\mapsto s+2\beta_{\cB}$ and $s\mapsto \frac{s}{1-2s\rho}$ coincide,
		and 
		\begin{equation}
			\label{e:deftaus}
			\tau(s):=\min\left\{s+2\beta_{\cB},\left(\frac{s}{1-2s\rho}\right)_{\!\!+}\right\}
			=\begin{cases}
				\frac{s}{1-2s\rho},&\text{if}\:\: s\le s^*;\\
				s+2\beta_{\cB},&\text{if}\:\: s> s^*,
			\end{cases}
		\end{equation}
		we have $r(s,\tau(s))=\min_{\tau\in\Sigma(s)}r(s,\tau)$.
		We consider the following two cases:
		
		$\bullet$ If $s\le s^*$, we deduce from \eqref{e:deftaus} that
		\begin{equation}
			r(s,\tau(s))=r\left(s,\frac{s}{1-2s\rho}\right)=
			1-2s\rho
		\end{equation}
		which is decreasing and, thus, optimized at $s=s^*$, obtaining
		\begin{equation}
			\label{e:case1}
			\widehat{r}=r(s^*,\tau(s^*))=1-2s^*\rho=\frac{s^*}{s^*+2\beta_{\cB}}=\frac{\sqrt{1+\frac{1}{\beta_{\cB}\rho}}-1}{\sqrt{1+\frac{1}{\beta_{\cB}\rho}}+1}
		\end{equation}
		with optimal $\widehat{\tau}:=\tau(s^*)$ and $\widehat{\theta}=s^*/\widehat{\tau}$.
		
		$\bullet$ If $s>s^*$, then \eqref{e:deftaus} yields
		\begin{equation}
			\label{e:rsqrtphi}
			r(s,\tau(s))=r\left(s,s+2\beta_{\cB}\right)=\sqrt{\phi(s)}, \textnormal{ where } \phi\colon s\mapsto \frac{(s-\beta_{\cA})^2+\beta_{\cA}(\frac{1}{\rho}-\beta_{\cA})}{(s+\beta_{\cA}+2\beta_{\cB})^2+\beta_{\cA}(\frac{1}{\rho}-\beta_{\cA})}.
		\end{equation}
		Note that, after simple computations, for every $s>0$,
		\begin{equation}
			\phi'(s)\ge 0\quad\Leftrightarrow\quad 
			s\ge\sqrt{\beta_{\cB}^2+2\beta_{\cA}\beta_{\cB}+\frac{\beta_{\cA}}{\rho}}-\beta_{\cB}=:s_1\ge 0.
		\end{equation}
		Therefore, if 
		$s^*\ge s_1$, $\phi$ is increasing on $[s^*,+\infty[$ and, thus, $\widehat{r}=r(s^*,\tau(s^*))$ is as in \eqref{e:case1}, in view of \eqref{e:rsqrtphi}.
		Otherwise, if $s^*< s_1$, 
		$\phi$ is decreasing on $[s^*,s_1]$ and increasing on $[s_1,+\infty[$, which implies 
		{\small\begin{equation}
				\label{e:ratecas2}
				\widehat{r}
				=r(s_1,\tau(s_1))
				=\frac{\sqrt{\beta_{\cA}\big(\frac{1}{\rho}-\beta_{\cA}\big)}}{\sqrt{\beta_{\cB}^2+2\beta_{\cA}\beta_{\cB}+
						\frac{\beta_{\cA}}{\rho}}+\beta_{\cA}+\beta_{\cB}}.
		\end{equation}}
		Note that $s^{*}=\sqrt{\beta_{\cB}^{2}+\frac{\beta_{\cB}}{\rho}}-\beta_{\cB}$ is exactly the value of $s_{1}$ obtained by
		replacing $\beta_{\cA}$ with $\beta_{\cB}/(1+2\beta_{\cB}\rho)$, and that $\beta_{\cA}\mapsto s_{1}(\beta_{\cA})$ is increasing. Hence
		$s^{*}\ge s_{1}\Leftrightarrow \beta_{\cA}\le\beta_{\cB}/(1+2\beta_{\cB}\rho)$, and the maximum of $s^{*}$ and $s_{1}$
		is obtained by evaluating $s_{1}$ at $\widehat{\beta}$ defined in \eqref{e:hatalpha}.
		Hence, both cases can be unified by taking $\widehat{s}=\max\{s^*,s_1\}=\sqrt{\beta_{\cB}^2+2\widehat{\beta}\beta_{\cB}+\frac{\widehat{\beta}}{\rho}}-\beta_{\cB}$ and it follows from \eqref{e:deftaus} that $\widehat{\tau}=\tau(\widehat{s})=2\beta_{\cB}+\widehat{s}$
		and, since $\widehat{\theta}=\widehat{s}/\widehat{\tau}$, we conclude \eqref{e:argmin}. The optimal rate can also be written in the unified form \eqref{e:main} by using $\widehat{\beta}$. Indeed, the result is clear from \eqref{e:ratecas2} if $\widehat{\beta}=\beta_{\cA}$. On the other hand, if $\widehat{\beta}=\beta_{\cB}/(1+2\beta_{\cB}\rho)$, by using $2u^{2}+3u+1=(2u+1)(u+1)$ at $u=\beta_{\cB}\rho$, we obtain 
		$$\sqrt{\beta_{\cB}^{2}+2\widehat{\beta}\beta_{\cB}+\frac{\widehat{\beta}}{\rho}}
		=\sqrt{\frac{\beta_{\cB}\left(2\beta_{\cB}^{2}\rho^{2}+3\beta_{\cB}\rho+1\right)}{\rho(1+2\beta_{\cB}\rho)}}
		=\sqrt{\frac{\beta_{\cB}(1+\beta_{\cB}\rho)}{\rho}}
		=s^*+\beta_{\cB}.$$ 
		Moreover, in this case $\widehat{\beta}\big(\frac{1}{\rho}-\widehat{\beta}\big)
		=\frac{\beta_{\cB}^{2}+\frac{\beta_{\cB}}{\rho}}{(1+2\beta_{\cB}\rho)^{2}}$, which yields
		\[
		\frac{\sqrt{\widehat{\beta}\big(\frac{1}{\rho}-\widehat{\beta}\big)}}
		{\sqrt{\beta_{\cB}^{2}+2\widehat{\beta}\beta_{\cB}+\frac{\widehat{\beta}}{\rho}}+\widehat{\beta}+\beta_{\cB}}
		=\frac{s^{*}+\beta_{\cB}}{s^{*}+3\beta_{\cB}+2\beta_{\cB}\rho(s^{*}+2\beta_{\cB})},
		\]
		and noting that 
		{$$s^*+2\beta_{\cB}=\sqrt{\beta_{\cB}^2+\frac{\beta_{\cB}}{\rho}}+\beta_{\cB}=\frac{\beta_{\cB}/\rho}{\sqrt{\beta_{\cB}^2+\frac{\beta_{\cB}}{\rho}}-\beta_{\cB}}=\frac{\beta_{\cB}}{\rho s^*}$$}
		we deduce
		\[
		\frac{\sqrt{\widehat{\beta}\big(\frac{1}{\rho}-\widehat{\beta}\big)}}
		{\sqrt{\beta_{\cB}^{2}+2\widehat{\beta}\beta_{\cB}+
				\frac{\widehat{\beta}}{\rho}}+\widehat{\beta}+\beta_{\cB}}
		=\frac{s^*(s^{*}+\beta_{\cB})}{(s^{*})^2+3\beta_{\cB} s^*+2\beta_{\cB}^2}=\frac{s^*}{s^{*}+2\beta_{\cB}},
		\]
		which coincides with \eqref{e:case1}.
		
		\ref{t:rrcontractive3}: 
		Since $\beta_{\cA}\rho\le1$, $\xi$ is well defined. For every $z=(x,y)\in\HH$,
		\begin{equation}\label{eq:exactconstants}
			\scal{\mathcal{A}z}{z}=\eta\|z\|^2
			\quad\text{and}\quad
			\|\mathcal{A}z\|^2=(\eta^2+\xi^2)\|z\|^2.
		\end{equation}
		Now, according to \eqref{eq:defmu} we have $\eta\ge\rho$ and 
		$\beta_{\cA}(\eta^2+\xi^2)\le\eta$. Indeed, if $\theta<(1+2\tau\rho)^{-1}$, then $\eta=\rho$, 
		$\xi=0$, and $\beta_{\cA}(\eta^2+\xi^2)=\beta_{\cA}\rho^2\le\rho=\eta$. If 
		$(1+2\tau\rho)^{-1}\le\theta\le1+2\beta_{\cA}/\tau$, then $\eta=\rho$ and 
		$\beta_{\cA}(\eta^2+\xi^2)=\rho$. If $\theta>1+2\beta_{\cA}/\tau$, then 
		$\eta=\beta_{\cA}^{-1}\ge\rho$, $\xi=0$, and $\beta_{\cA}(\eta^2+\xi^2)=\eta$. Hence, since $\mathcal{A}$ is linear, it follows from \eqref{eq:exactconstants} that it is $\rho$-strongly monotone and 
		$\beta_{\cA}$-cocoercive. Moreover, since $\mathcal{A}$ is monotone and continuous it is maximally monotone.  On the other hand,
		$\mathcal{B}=0$ is $\beta_{\cB}$-cocoercive for every $\beta_{\cB}\in\RPP$. We have 
		$Z=\mathcal{A}^{-1}(\{0\})=\{(0,0)\}$.
		
		Set now $\alpha=(2+2\beta_{\cA}/\tau)^{-1}$, $\beta=1+\tau\rho$, and $\lambda=1+\theta$. 
		As in the proof of \ref{t:rrcontractive1}, applied to $\tau\mathcal{A}$, we have 
		$2\beta/(2\beta-1)=1+(1+2\tau\rho)^{-1}$ and $1/\alpha=1+(1+2\beta_{\cA}/\tau)$. 
		Furthermore, since 
		$(\Id+\tau\mathcal{A})(x,y)=\big((1+\tau\eta)x-\tau\xi y,(1+\tau\eta)y+\tau\xi x\big)$, 
		we obtain
		\begin{equation*}
			J_{\tau\mathcal{A}}\colon(x,y)\mapsto\frac{1}{\delta}
			\begin{pmatrix}
				(1+\tau\eta)x+\tau\xi y\\
				(1+\tau\eta)y-\tau\xi x
			\end{pmatrix},
			\quad\text{where}\quad
			\delta=(1+\tau\eta)^2+\tau^2\xi^2,
		\end{equation*}
		which is the operator in \eqref{e:tightT} with $\mu=(1+\tau\eta)/\delta$ and 
		$\sigma=\tau\xi/\delta$. We claim that these values of $\mu$ and $\sigma$ match the expressions in \eqref{e:tightmusigma}. Indeed, if $\theta<(1+2\tau\rho)^{-1}$, then $\mu=(1+\tau\rho)^{-1}=\beta^{-1}$ and $\sigma=0$. 
		If $\theta>1+2\beta_{\cA}/\tau$, then $\mu=\beta_{\cA}/(\tau+\beta_{\cA})=1-2\alpha$ 
		and $\sigma=0$. Otherwise, we have $1+\tau\eta=\beta$ and
		\begin{equation}\label{e:deltaproof}
			\delta=\beta^2+\tau^2\xi^2=1+2\tau\rho+\frac{\tau^2\rho}{\beta_{\cA}}
			=\frac{2\beta(1-\alpha)-1}{1-2\alpha}.
		\end{equation}
		Then, it follows from \eqref{e:deltaproof} that 
		$\sigma^2=(\delta-\beta^2)/\delta^2=\mu/\beta-\mu^2$ and that $\mu=\beta/\delta=\beta(1-2\alpha)/(2\beta(1-\alpha)-1)$. 
		Altogether, $R_{\tau\mathcal{A}}^{\theta}=(1+\theta)J_{\tau\mathcal{A}}-\theta\Id=S_{\lambda}$, 
		and Proposition~\ref{p:reflectcontr}.\ref{p:reflectcontr2} yields 
		$\|R_{\tau\mathcal{A}}^{\theta}z\|=\rate_T(\lambda)\|z\|=\rate_{\tau\mathcal{A}}(\theta)\|z\|$ 
		for every $z\in\HH$. Finally, if $\theta\in I(\tau)$, then $J_{\theta\tau\mathcal{B}}=\Id$, 
		so $R_{\theta\tau\mathcal{B}}^{1/\theta}=\Id$ and 
		$R_{\theta\tau\mathcal{B}}^{1/\theta}\circ R_{\tau\mathcal{A}}^{\theta}=R_{\tau\mathcal{A}}^{\theta}$, 
		which attains $\rate(\theta,\tau)=\rate_{\tau\mathcal{A}}(\theta)$.
	\end{proof}
	
	The following result is a consequence of 
	Proposition~\ref{t:genDRtau}.\ref{t:genDRtau1} and Theorem~\ref{t:rrcontractive}.\ref{t:rrcontractive2}.
	
	\begin{theorem}\label{t:strongconvergenceDR}
		In the context of Problem~\ref{prob:main}, assume 
		\eqref{e:strongmon}. Let $\tau \in \RPP$, let $\theta \in \RPP$, let $z_0 \in \HH$, and consider the recurrence
		\begin{align}
			(\forall n\in\NN)\hspace{0.2cm}
			&\left\lfloor 
			\begin{array}{ll}
				x_{n}=J_{\tau \mathcal{A}}z_n\\
				y_n=J_{\theta\tau\mathcal{B}}\big((1+\theta)x_{n}
				-\theta z_n\big)\\
				z_{n+1}=z_n+\left(\dfrac{1+\theta}{\theta}\right)
				(y_n-x_{n}).
				\label{e:RDR2}
			\end{array}
			\right.
		\end{align}
		Then, if $\beta_{\cA}+\beta_{\cB}>0$ and $\theta \in I(\tau)$, where $I(\tau)$ is defined in \eqref{e:I(tau)}, the sequences $(z_{n})_{n \in \NN}$ and $(x_{n})_{n \in \NN}$ converge strongly to the unique fixed point of $R_{\theta \tau\mathcal{B}}^{1/\theta}\circ 
		R_{\tau\mathcal{A}}^{\theta}$ and the unique point in $Z$, respectively, with a linear convergence rate given by $\rate(\theta,\tau)$ defined in \eqref{e:defomega}. 
		Moreover, this rate is optimized for $(\widehat{\theta},\widehat{\tau})$ in \eqref{e:argmin}, which leads to the optimal rate in \eqref{e:main}.
	\end{theorem}
	\begin{proof} Note that \eqref{e:RDR2} is a particular instance of 
		\eqref{e:RDR} when $\lambda_n \equiv 1$. Then, as in the proof of 
		Proposition~\ref{t:RDRweakc}, for every $n \in \NN$, $z_{n+1} = 
		R_{\theta \tau\mathcal{B}}^{1/\theta}\circ 
		R_{\tau\mathcal{A}}^{\theta} z_n$. It follows from Theorem~\ref{t:rrcontractive}.\ref{t:rrcontractive2} that $R_{\theta \tau\mathcal{B}}^{1/\theta}\circ 
		R_{\tau\mathcal{A}}^{\theta}$ is a contraction with Lipschitz constant given by $\rate(\theta,\tau)$ defined in \eqref{e:defomega}. Therefore $\Fix (R_{\theta \tau\mathcal{B}}^{1/\theta}\circ 
		R_{\tau\mathcal{A}}^{\theta})$ is a singleton, say $\{z\}$,
		in view of \cite[Theorem~1.50]{bauschke2011:Convex_analysis}, and, for every $n\in\NN$,
		$\|z_n-z\|\le \omega(\theta,\tau)^n\|z_0-z\|$. Finally, since Proposition~\ref{t:genDRtau}.\ref{t:genDRtau1} asserts that $Z = \{J_{\tau\mathcal{A}}z\}$ and 
		$J_{\tau\mathcal{A}}$ is nonexpansive, 
		\begin{equation*}
			\|x_n - J_{\tau\mathcal{A}}z\|= \|J_{\tau\mathcal{A}}z_n-J_{\tau\mathcal{A}}z\|\leq \|z_n-z\|
		\end{equation*}
		and the result follows by Theorem~\ref{t:rrcontractive}.\ref{t:rrcontractive2}.
	\end{proof}
	\begin{remark} 
		\begin{enumerate}
			\item In the context of Theorem~\ref{t:strongconvergenceDR}, set $\theta=1$. Then $T=R_{\tau\mathcal{B}}^{1}\circ R_{\tau\mathcal{A}}^{1}$
			is the composition of the exact reflections,
			\eqref{e:RDR2} corresponds to the
			Peaceman--Rachford algorithm, and we recover the rate in \cite{Giselsson2017} with optimal relaxation step. Indeed, when $\rho>0$ and $\beta_{\cA}>0$, 
			from \eqref{e:defomega} we
			obtain 
			\vspace{-.4cm}
			{\small\begin{equation*}
					\rate(1,\tau)=
					\sqrt{1-\dfrac{4\tau\beta_{\cA}\rho}{\beta_{\cA}+\tau^2\rho+2\tau\beta_{\cA}\rho}},
			\end{equation*}}
			which coincides with the constant obtained in 
			\cite[Theorem~7.4]{Giselsson2017} in the case when $\beta_{\cB}=0$. It is optimized at $\tau^*=\sqrt{\beta_{\cA}/\rho}$ with optimal value
			{\small\begin{equation}
					\label{e:optgiss}
					\rate(1,\tau^*)=\sqrt{\frac{1-\sqrt{\beta_{\cA}\rho}}{1+\sqrt{\beta_{\cA}\rho}}}.
			\end{equation}}
			In general, the rate obtained in Theorem~\ref{t:rrcontractive} for optimal parameters is strictly lower than the best possible rate for Douglas--Rachford when $\beta_{\cB}>0$. Indeed, $\rate(\widehat{\theta},\widehat{\tau})$ in \eqref{e:main} is strictly decreasing in $\beta_{\cB}$ and coincides with \eqref{e:optgiss} when $\beta_{\cB}=0$. 
			Figure~\ref{fig:diff} shows a color map of the difference between 
			$\rate(\widehat{\theta},\widehat{\tau})$ and 
			$\rate(1,\tau^*)$ for different values of $\rho\in\left[0,1\right]$ and $\beta_{\cB}>0$, where we set $\beta_{\cA}=1$, without loss of generality.

			\item 
			When $\beta_{\cB}>0$ and $\beta_{\cA}=0$, \cite[Theorem~5.6]{Giselsson2017}
			provides a tight worst-case linear convergence rate for Douglas--Rachford, given by $\frac{1}{1+\sqrt{\beta_{\cB}\rho}}$ for a specific relaxation step and step-size. This bound is strictly larger than the rate provided in Theorem~\ref{t:rrcontractive} in this case, more precisely,
			\begin{equation}
				\frac{1}{1+\sqrt{\beta_{\cB}\rho}}>\frac{\sqrt{1+\beta_{\cB}\rho}-\sqrt{\beta_{\cB}\rho}}{\sqrt{1+\beta_{\cB}\rho}+\sqrt{\beta_{\cB}\rho}}.
			\end{equation}
			Moreover, the composition of relaxed resolvents is a strict contraction, and the relaxation does not improve the rate. In Figure~\ref{fig:comp1} we compare 
			both rates when, without loss of generality, we set $\beta_{\cB}=1$ and $\rho>0$.
			\begin{figure}[htbp]
				\centering
				\begin{subfigure}[b]{0.48\textwidth}
					\centering
					\includegraphics[scale=0.25]{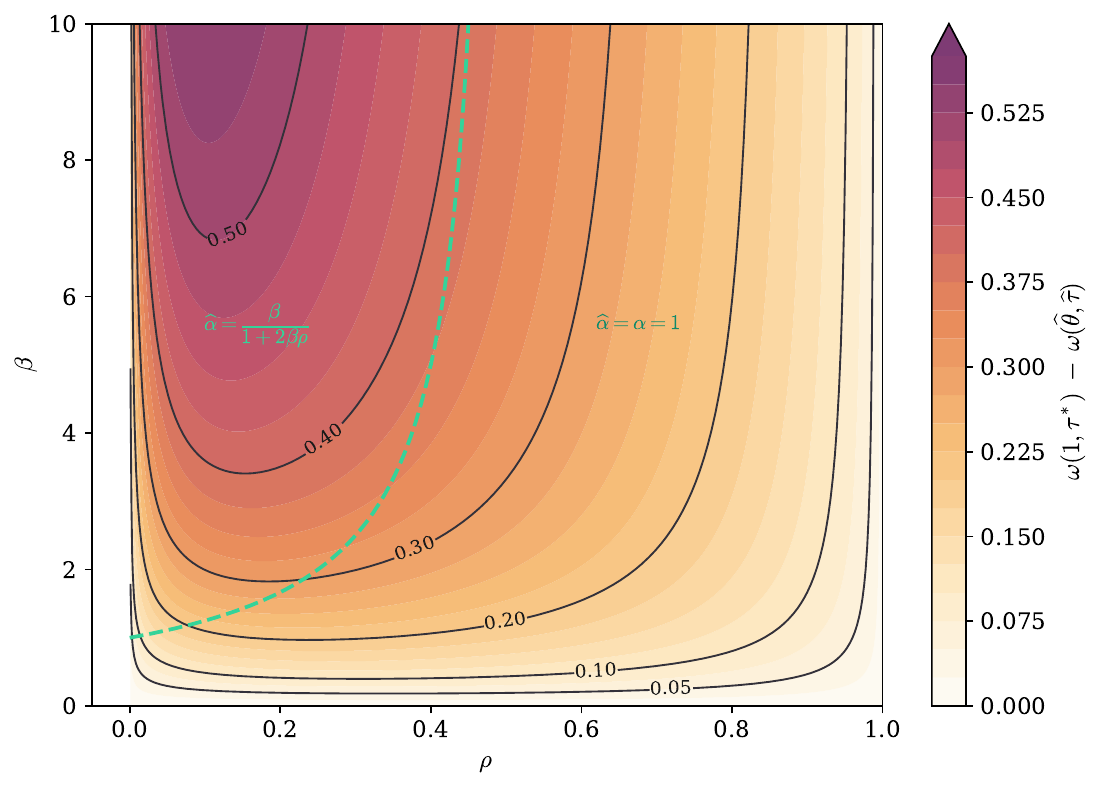}
					\caption{Case $\beta_{\cB}>0$, $\beta_{\cA}=1$, and $\rho\in\left[0,1\right]$.}
					\label{fig:diff}
				\end{subfigure}
				\hfill 
				\begin{subfigure}[b]{0.48\textwidth}
					\centering
					\includegraphics[scale=0.25]{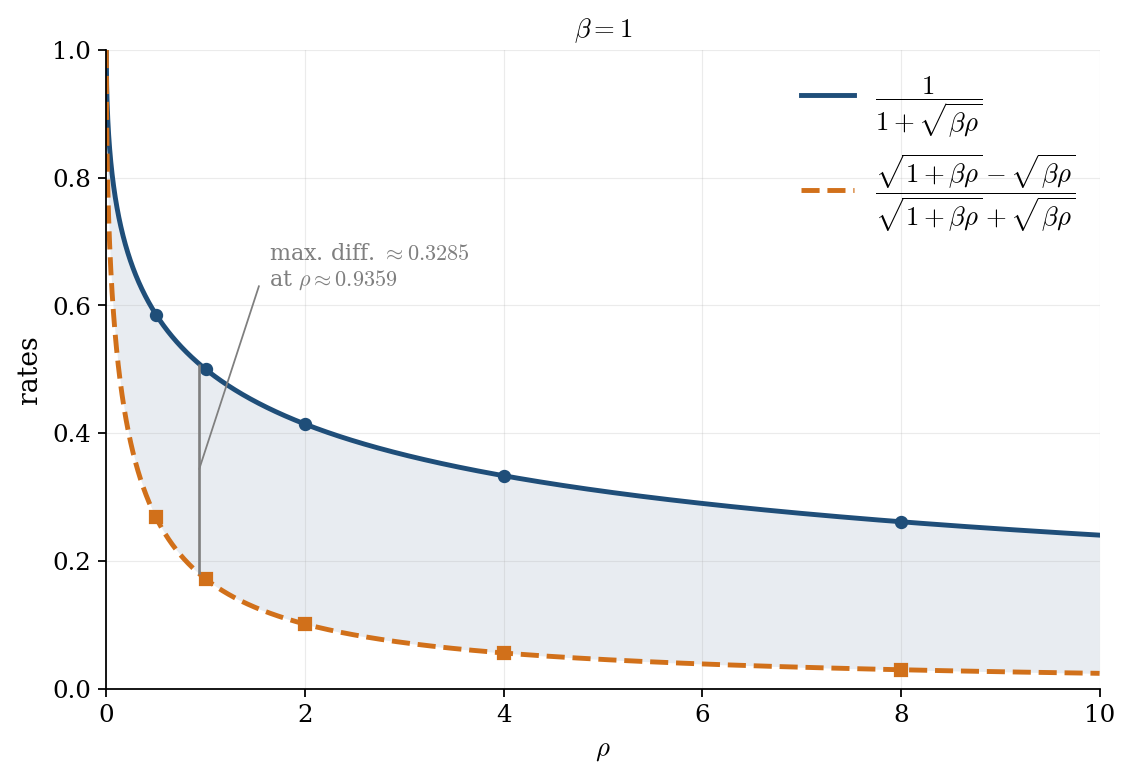}
					\caption{Case $\beta_{\cB}=1$ and $\beta_{\cA}=0$.}
					\label{fig:comp1}
				\end{subfigure}
				\caption{Comparison of optimal rates in \cite{Giselsson2017} and Theorem~\ref{t:rrcontractive}: (a) using \cite[Theorem~7.4]{Giselsson2017}; (b) using \cite[Theorem~5.6]{Giselsson2017}.}
				\label{fig:combined_comparison}
				\vspace{-.4cm}
			\end{figure}
			
			\item In the classical case $\theta=1$, the weak convergence of the shadow
			sequence $(x_n)_{n\in\NN}$ was established in \cite{Svai11}.
			Proposition~\ref{t:RDRweakc} extends this property to every 
			admissible $\theta$.
			
			\item 
			Suppose that \eqref{e:strongmon} holds, that $\beta_{\cA}>0$ and that
			$\beta_{\cB}=0$. Then $\mathcal{A}$ is $\rho$-strongly monotone and
			$\beta_{\cA}^{-1}$-Lipschitz continuous, so that
			\cite[Theorem~4.8(i)]{DaoPhan2019} applies with monotonicity constants
			$\rho$ and $0$ and yields the linear convergence of the relaxed iteration
			$z_{n+1}=(1-\kappa)z_n+\kappa\,R_{\theta\tau\mathcal{B}}^{1/\theta}
			(R_{\tau\mathcal{A}}^{\theta}z_n)$ with $\kappa\in\left]0,1\right[$ and
			$\theta\in\left[(1+2\tau\rho)^{-1},1\right]$. In addition, the lower endpoint of this
			interval coincides with the threshold $(1+2\tau\rho)^{-1}$ in \eqref{e:defomega}.
			Theorem~\ref{t:rrcontractive}.\ref{t:rrcontractive2} differs in three
			respects: it concerns the unrelaxed composition, it allows $\theta$ up to
			$1+2\beta_{\cA}/\tau$ by exploiting the cocoercivity of $\mathcal{A}$, and it
			provides the explicit constant $\rate(\theta,\tau)$ together with its
			minimizers \eqref{e:argmin}. Since the constant of
			\cite[Theorem~4.8]{DaoPhan2019} is obtained under Lipschitz continuity
			instead of cocoercivity and concerns a relaxed operator with $\kappa<1$, the
			two constants are not directly comparable. However, for $\theta=1$, the constant $\rate(1,\tau)$ is smaller than the one obtained
			by applying \cite[Corollary~4.9]{DaoPhan2019} with Lipschitz constant
			$\beta_{\cA}^{-1}$. It is strictly smaller whenever $\beta_{\cA}\rho<1$. This is consistent with
			the fact that cocoercivity is a more restrictive assumption than Lipschitz continuity.
			
			\item Tight contraction factors for the classical Douglas--Rachford
			algorithm, together with optimal parameter choices, are obtained via
			performance estimation in \cite{RyuTaylorBergelingGiselsson2020} for
			several combinations of strong monotonicity, cocoercivity, and Lipschitz
			continuity, including the case in which $\mathcal{A}$ is strongly monotone
			and $\mathcal{B}$ is cocoercive. The case in which one operator is
			Lipschitz continuous and the other is strongly monotone is studied in
			\cite{MoursiVandenberghe2019}. Since these results concern relaxed
			iterations and, in the latter case, a different regularity class, the
			corresponding constants are not directly comparable with
			$\rate(\theta,\tau)$.
		\end{enumerate}
		
	\end{remark}

\begin{figure}[htbp]
			\centering
			\begin{subfigure}[b]{0.32\textwidth}
				\centering
				\includegraphics[width=\textwidth]{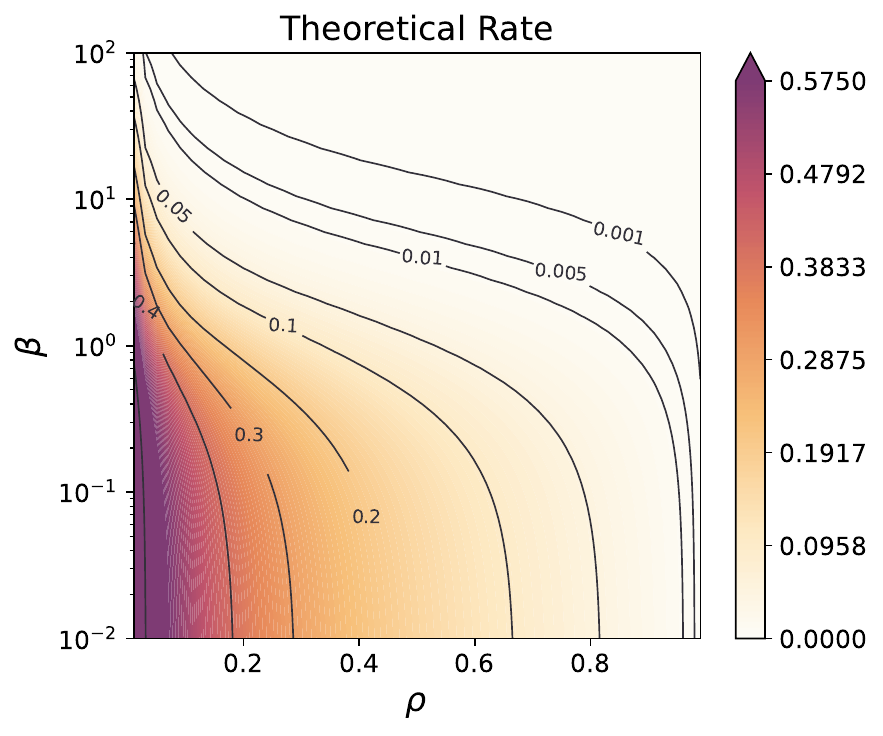}
			\end{subfigure}
			\hfill
			\begin{subfigure}[b]{0.32\textwidth}
				\centering
				\includegraphics[width=\textwidth]{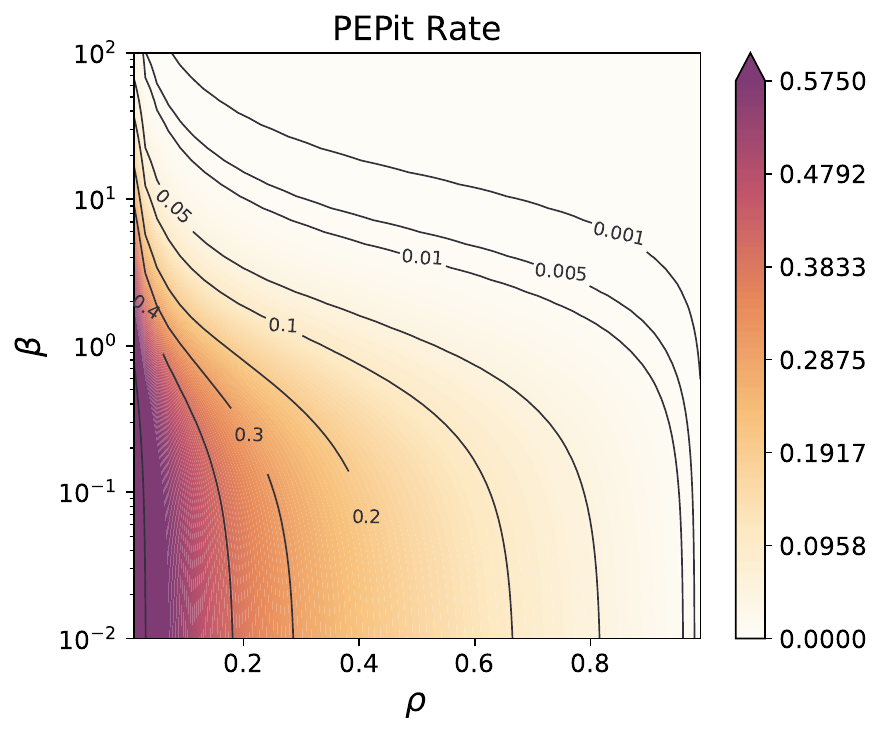}
			\end{subfigure}
			\hfill
			\begin{subfigure}[b]{0.32\textwidth}
				\centering
				\includegraphics[width=\textwidth]{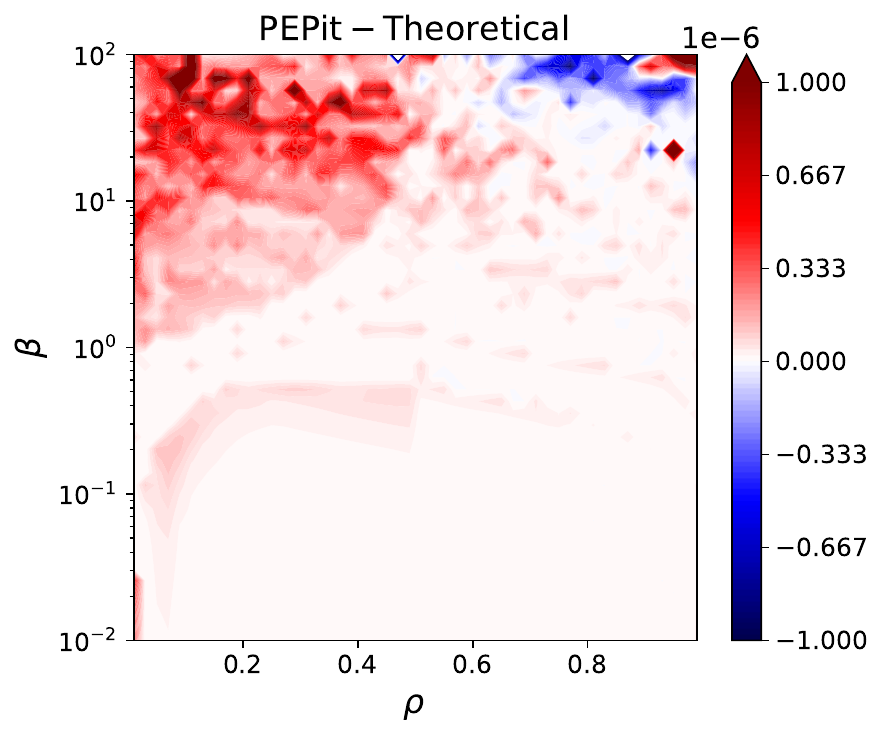}
			\end{subfigure}
			\caption{Comparison between the theoretical convergence rate and the numerical worst-case obtained via PEPit for the Relaxed Douglas-Rachford algorithm as a function of the parameters $\rho\in\left[0,\beta_{\cA}\right]$ and $\beta_{\cB}>0$ (log scale). Without loss of generality, we set $\beta_{\cA} = 1$.}
			\label{fig:rdr_comparison}
		\end{figure}

		We conclude this section by comparing our theoretical rate with 
		the computer-assisted worst-case analysis of PEPit \cite{GoujaudEtAl2024PEPit}, which finds its roots in the theoretical results of \cite{DroriTeboulle2014} (see also \cite{TaylorHendrickxGlineur2017,TaylorHendrickxGlineur2017Smooth}).
		In Figure~\ref{fig:rdr_comparison} we observe a difference of less 
		than $10^{-6}$ between our theoretical rate and that obtained by 
		PEPit.
		
		\section{Convex optimization}\label{sec:convex}
		Consider the following problem.
		\begin{problem}
			\label{prob:opti}
			Let $\beta_{f}\ge 0$, let $\beta_{g}\ge 0$, let $\rho\ge 0$, and let 
			$f$ and $g$ be two functions in $\Gamma_0(\HH)$ such that 
			$f$ is $\rho$-strongly convex, $\partial f$ is 
			$\beta_{f}$-cocoercive, and $\partial g$ is 
			$\beta_{g}$-cocoercive. The problem is 
			\begin{equation}
				\min_{x\in\HH}f(x)+g(x),
			\end{equation}
			under the assumption that solutions exist.
		\end{problem}
		
		Note that, in the case when $\beta_{f}>0$ and $\dom\partial f=\HH$, the $\beta_{f}$-cocoercivity of $\partial f$ implies the G\^ateaux differentiability of $f$ \cite[Propositions~16.27 \& 17.31]{bauschke2011:Convex_analysis}.	In addition, it follows from the Cauchy--Schwarz inequality, that $\nabla f$ is $\beta_{f}^{-1}$-Lipschitz continuous, which implies the Fr\'echet differentiability of $f$ \cite[Corollary~17.42]{bauschke2011:Convex_analysis}. Moreover, by the Baillon-Haddad theorem,  if $f$ is Fréchet differentiable, $\nabla f$ is $\beta_{f}$-cocoercive if and 
		only if $\nabla f$  is $\beta_{f}^{-1}$-Lipschitz continuous. An 
		analogous equivalence for $\nabla g$ is obtained when $\beta_{g}>0$ and  $\dom\partial g=\HH$. These additional properties in the context of convex functions can be exploited 
		to obtain stronger contractive properties of the relaxed 
		resolvent of $\partial f$. Since $J_{\partial f} = \prox_{f}$, this operator reduces to
		\begin{equation}
			(\forall \theta\ge -1)\quad R_{\partial f}^{\theta}=(1+\theta)\prox_f-\theta\Id
		\end{equation}
		which is simply denoted by $R_{f}^{\theta}$. The next result provides an improved Lipschitz constant for $R_{f}^{\theta}$ under the strong convexity of $f$.
		
		\begin{theorem}
			\label{p:optiref}
			In the context of Problem~\ref{prob:opti}, the following 
			holds:
			\begin{enumerate}
				
				\item\label{p:optirefi}
				For every $\theta\ge -1$, 
				$R_{f}^{\theta}$ is 
				$\widehat{\rate}_f(\theta)$-Lipschitz continuous, where
				\begin{equation}\label{eq:omega_f}
					\widehat{\rate}_f(\theta)=
					\max\left\{\frac{1-\theta\rho}{1+\rho},
					\frac{\theta-\beta_{f}}{1+\beta_{f}}\right\}.
				\end{equation}
				

				\item\label{p:optirefii}
				Let $\tau>0$ and define
				\begin{equation}\label{e:I(tau)tilde}
					\widehat{I}(\tau)=\left]0,1+\dfrac{2\beta_{f}}{\tau}\right]\cap\left[1-\dfrac{2\beta_{g}}{\tau},+\infty\right[.
				\end{equation}
				For every $\theta\in \widehat{I}(\tau)$,
				$R_{\theta \tau g}^{1/\theta}\circ 
				R_{\tau f}^{\theta}$ is
				$\widehat{\rate}(\theta,\tau)$-Lipschitz 
				continuous, 
				where 
				\begin{equation}
					\label{e:r1}
					\widehat{\rate}(\theta,\tau)=\max\left\{
					\frac{1-\tau\theta\rho}{1+\tau\rho},\frac{\tau\theta-\beta_{f}}{\tau+\beta_{f}}\right\}
					\in\left]0,1\right].
				\end{equation}
				
				Moreover, by 
				we have
				\begin{equation}\label{e:mainfg}
					\widehat{\rate}(\theta^*,\tau^*)=
					\min_{\tau\in\left]0,+\infty\right[}\;
					\min_{\theta\in \widehat{I}(\tau)}\;\widehat{\rate}(\theta,\tau)
					=\frac{\sqrt{1+\beta_{g}\rho}-\sqrt{\rho(\beta_{f}+\beta_{g})}}
					{\sqrt{1+\beta_{g}\rho}+\sqrt{\rho(\beta_{f}+\beta_{g})}},
				\end{equation}
				where 
				\begin{equation}\label{e:argminfg}
					\tau^*=\beta_{g}+\sqrt{\frac{(\beta_{f}+\beta_{g})(1+\beta_{g}\rho)}{\rho}}
					\quad\text{and}\quad
					\theta^*=\frac{\sqrt{(\beta_{f}+\beta_{g})(1+\beta_{g}\rho)}-\beta_{g}\sqrt{\rho}}
					{\sqrt{(\beta_{f}+\beta_{g})(1+\beta_{g}\rho)}+\beta_{g}\sqrt{\rho}}.
				\end{equation}
				
				\item\label{p:optirefiii}Let $\tau>0$ and $\theta\in\widehat{I}(\tau)$, and suppose that 
				$\beta_{f}>0$. The constant
				$\widehat{\rate}(\theta,\tau)$ in \eqref{e:r1} is
				attained by taking 
				$\HH=H_1\oplus H_2\oplus H_3$ and
				$\widehat{f}\colon 
				(x_1,x_2,x_3)\mapsto\frac{\rho}{2}(\|x_1\|^2+\|x_2\|^2)
				+\frac{1}{2\beta_{f}}\|x_3\|^2$, where $H_1$, $H_2,$ and $H_3$ are real Hilbert spaces. In addition, if we assume $\beta_g>0$,
				the optimal Lipschitz constant $\widehat{\omega}(\theta^*,\tau^*)$ in \eqref{e:mainfg} is attained by taking 
				$\widehat{g}\colon (x_1,x_2,x_3)\mapsto \frac{1}{2\beta_g}\|x_1\|^2$, and, if we assume $\beta_g=0$, by taking $\widehat{g}\equiv0$.
			\end{enumerate}
		\end{theorem}
		\begin{proof}
			\ref{p:optirefi}: Set 
			$\phi=f+\frac{1}{2}\|\cdot\|^2$ 
			and note that 
			$\beta_{f}/(1+\beta_{f})\le 1/(1+\rho)$ since $\beta_{f}\rho\le 1$. It follows from 
			\cite[Proposition~12.30, Proposition~14.1, 
			Remark~14.4]{bauschke2011:Convex_analysis} that $\phi^*=f^*\,\square\, (\frac{1}{2}\|\cdot\|^2)$ 
			is real-valued and Fr\'echet differentiable, with $\nabla\phi^*=\prox_f$. 
			Since $\phi^{**}=\phi$ \cite[Corollary~13.38]{bauschke2011:Convex_analysis} 
			and $\phi-\frac{1}{2}(1+\rho)\|\cdot\|^2=f-\frac{\rho}{2}\|\cdot\|^2$ is convex, 
			\cite[Proposition~14.2]{bauschke2011:Convex_analysis} applied to $\phi^*$ 
			implies that $\frac{1}{2(1+\rho)}\|\cdot\|^2-\phi^*$ is convex. Let us show that $h=\phi^*-\frac{\beta_f}{2(1+\beta_f)}\|\cdot\|^2$ is 
			convex. This is clear if $\beta_{f}=0$. If $\beta_{f}>0$, then 
			$\partial f$ is single-valued on its domain and 
			\cite[Example~22.7 and Proposition~20.22]{bauschke2011:Convex_analysis} 
			imply that $(\partial f)^{-1}$ is maximally monotone and 
			$\beta_{f}$-strongly monotone. Hence, 
			\cite[Proposition~22.11(ii)]{bauschke2011:Convex_analysis} yields 
			$\dom\partial f=\HH$ and $f$ 
			is Fr\'echet differentiable with a $\beta_{f}^{-1}$-Lipschitz continuous 
			gradient. Therefore, \cite[Theorem~18.15]{bauschke2011:Convex_analysis} 
			implies that $\frac{1+\beta_{f}}{2\beta_f}\|\cdot\|^2-\phi=\frac{1}{2\beta_f}\|\cdot\|^2-f$ is convex, and 
			\cite[Proposition~14.2]{bauschke2011:Convex_analysis} applied to $\phi$ 
			yields the convexity of $h$. Altogether, $h$ and $(1/(1+\rho)-\beta_{f}/(1+\beta_{f}))\|\cdot\|^2/2-h=\|\cdot\|^2/(2(1+\rho))-\phi^*$ are 
			convex. If $\beta_{f}/(1+\beta_{f})<1/(1+\rho)$, \cite[Theorem~18.15]{bauschke2011:Convex_analysis} 
			implies that $\nabla h=\prox_f-\frac{\beta_{f}}{1+\beta_{f}}\Id$ is $\big(\frac{1}{1+\rho}-\frac{\beta_f}{1+\beta_f}\big)^{-1}$-cocoercive and, if 
			$\beta_{f}/(1+\beta_{f})=1/(1+\rho)$, $h$ is affine and $\nabla h$ is constant. In both cases, for every 
			$x$ and $y$ in $\HH$, 
			we obtain
			\begin{equation}
				\label{e:auxopt}
				\scalebox{0.95}{$\left\|\nabla \phi^*(x)-\nabla \phi^*(y)-\frac{\beta_f(x-y)}{1+\beta_f}\right\|^2\le\left(\frac{1}{1+\rho}-\frac{\beta_f}{1+\beta_f}\right)\scal{\nabla \phi^*(x)-\nabla \phi^*(y)-\frac{\beta_f(x-y)}{1+\beta_f}}{x-y}$},
			\end{equation}
			which yields 
			$$\scalebox{0.95}{$\|\nabla \phi^*(x)-\nabla \phi^*(y)\|^2\le\left(\frac{\beta_f}{1+\beta_f}+\frac{1}{(1+\rho)}\right)\scal{\nabla \phi^*(x)-\nabla \phi^*(y)}{x-y}-\frac{\beta_f}{(1+\beta_f)(1+\rho)}\|x-y\|^2$}$$
			and,  since $\frac{\scal{\nabla \phi^*(x)-\nabla \phi^*(y)}{x-y}}{\|x-y\|}\leq \|\nabla \phi^*(x)-\nabla \phi^*(y)\|$, we conclude
			\begin{equation}\label{eq:upperboundprox}
				\frac{\beta_f}{(1+\beta_f)}\|x-y\|^2\le\scal{\nabla \phi^*(x)-\nabla \phi^*(y)}{x-y}\le \frac{1}{1+\rho}\|x-y\|^2.
			\end{equation}
			Hence, since 
			$R_f^{\theta}x-R_f^{\theta}y=(1+\theta)(\nabla \phi^*(x)-\nabla \phi^*(y))-\theta(x-y)$, we have
			\begin{align*} 
				\scalebox{0.95}{$\|R_f^{\theta}x-R_f^{\theta}y\|^2
					\le(1+\theta)$}&\scalebox{0.95}{$\left((1+\theta)\left(\frac{\beta_f}{1+\beta_f}+\frac{1}{1+\rho}\right)-2\theta\right)\scal{\nabla \phi^*(x)-\nabla \phi^*(y)}{x-y}$}
				\\
				&\scalebox{0.95}{$+\left(\theta^2-(1+\theta)^2\frac{\beta_f}{(\beta_f+1)(1+\rho)}\right)\|x-y\|^2$}.
			\end{align*}
			Since $1+\theta\ge0$, if $(1+\theta)\left(\frac{\beta_f}{1+\beta_f}+\frac{1}{1+\rho}\right)-2\theta\ge0$, \eqref{eq:upperboundprox} yields 
			$$\scalebox{0.95}{$\|R_f^{\theta}x-R_f^{\theta}y\|\le\left|\frac{1+\theta}{1+\rho}-\theta\right|\,\|x-y\|
				=\frac{|1-\theta\rho|}{1+\rho}\|x-y\|$}.$$ Otherwise, the \eqref{eq:upperboundprox}
			yields 
			$$\scalebox{0.95}{$\|R_f^{\theta}x-R_f^{\theta}y\|\le\left|\frac{(1+\theta)\beta_f}{1+\beta_f}-\theta\right|\,\|x-y\|
				=\frac{|\theta-\beta_{f}|}{1+\beta_{f}}\|x-y\|$}.$$ Therefore, 
			by combining the last two cases,
			the Lipschitz constant in \eqref{e:r1} for $R_f^{\theta}$ is obtained by noting that 
			$\frac{1-\theta\rho}{1+\rho}+\frac{\theta-\beta_{f}}{1+\beta_{f}}
			=\frac{(1+\theta)(1-\beta_{f}\rho)}{(1+\rho)(1+\beta_{f})}\ge0$ 
			since
			$\beta_{f}\rho\le1$ and $\theta\ge-1$.
			\ref{p:optirefii}: Note that $\tau f$ is $\tau\rho$-strongly convex and $\partial (\tau f)$ is
			$(\beta_{f}/\tau)$-cocoercive, while
			$\theta\tau 
			\partial g$ 
			is  $(\beta_{g}/(\theta\tau))$-cocoercive. 
			Hence,
			it follows from \ref{p:optirefi}
			that, for every $\theta>0$, $R_{\tau f}^{\theta}$ is 
			$\max\left\{\frac{1-\theta\tau\rho}{1+\tau\rho},
			\frac{\tau\theta-\beta_{f}}{\tau+\beta_{f}}\right\}$-Lipschitz
			continuous and $R_{\theta\tau g}^{1/\theta}$ is 
			$\max\left\{1,
			\frac{\tau-\beta_{g}}{\theta\tau+\beta_{g}}\right\}$-Lipschitz
			continuous. Hence, $R_{\theta \tau g}^{1/\theta}\circ 
				R_{\tau f}^{\theta}$ is
				$\widehat{\rate}(\theta,\tau)$-Lipschitz 
			\begin{equation}
				\widehat{\rate}(\theta,\tau)=\max\left\{
				\frac{1-\theta\tau\rho}{1+\tau\rho},\frac{\tau\theta-\beta_{f}}{\tau+\beta_{f}}\right\}\max\left\{1,
				\frac{\tau-\beta_{g}}{\theta\tau+\beta_{g}}\right\}.
			\end{equation}
			Noting that 
			\begin{equation}
				\label{e:leq1}
				\frac{\tau-\beta_{g}}{\theta\tau+\beta_{g}}\le 
				1\quad\Leftrightarrow\quad(1-\theta)\tau\le 
				2\beta_{g}
			\end{equation}
			and 
			\begin{equation}\label{e:leq2}
				\max\left\{
				\frac{1-\theta\tau\rho}{1+\tau\rho}
				,\frac{\tau\theta-\beta_{f}}{\tau+\beta_{f}}\right\}\le 1
				\quad\Leftrightarrow\quad (\theta-1)\tau\le 2\beta_{f},
			\end{equation}
			we have 
			$$(1-\theta)\tau\in\left[-2\beta_{f}, 
			2\beta_{g}\right]\quad\Leftrightarrow\quad 
			\widehat{\rate}(\theta,\tau)=\max\left\{
			\frac{1-\theta\tau\rho}{1+\tau\rho},
			\frac{\tau\theta-\beta_{f}}{\tau+\beta_{f}}\right\}\in\left]0,1\right].$$
			Note that, since $\Theta_1 \colon \theta\mapsto\frac{1-\theta\tau\rho}{1+\tau\rho}$ and $\Theta_ 2 \colon \theta \mapsto \frac{\tau\theta-\beta_{f}}{\tau+\beta_{f}}$  are decreasing and increasing, respectively, $\theta \mapsto \max\{\Theta_1(\theta),\Theta_2(\theta)\}$ achieves its global minimum at $\bar\theta$ satisfying $\Theta_1(\bar\theta)=\Theta_2(\bar\theta)$, that is, $\bar\theta(\tau)=({1+\beta_{f}\rho+\frac{2\beta_{f}}{\tau}})/
			({1+\beta_{f}\rho+2\tau\rho})$. Hence, $$\min_{\theta \in 
				\RR}\widehat{\rate}(\theta,\tau)=\widehat{\rate}(\bar\theta(\tau),\tau)
			=\frac{\tau\bar\theta(\tau)-\beta_{f}}{\tau+\beta_{f}}=\frac{1-\beta_{f}\rho}
			{1+\beta_{f}\rho+2\tau\rho}.$$
			Since $\tau \mapsto \frac{1-\beta_{f}\rho}
			{1+\beta_{f}\rho+2\tau\rho}$ is decreasing and \eqref{e:leq2} holds directly, $\min_{\tau>0}\min_{\theta\in 
				\widehat{I}(\tau)}\widehat\rate(\theta,\tau)$ is achieved for the largest $\tau>0$ 
			guaranteeing \eqref{e:leq1}. More precisely,
			$$(1-\theta(\tau))\tau\leq 2\beta_{g}\quad\Leftrightarrow\quad 
			\frac{\tau^2\rho-\beta_{f}}{1+\beta_{f}\rho+2\tau\rho}\le 
			\beta_{g}\quad\Leftrightarrow
			\quad 
			\tau\in\left]0,\tau^*\right],$$
			where 
			$\tau^*
			$ is defined in \eqref{e:argminfg}. Hence, $\min_{\tau>0}\min_{\theta\in 
				I(\tau)}\rate(\theta,\tau)=\rate(\theta(\tau^*),\tau^*)$ and
			$$
			\rate(\theta(\tau^*),\tau^*)
			=\frac{1-\beta_{f}\rho}{(\sqrt{1+\beta_{g}\rho}+\sqrt{\rho(\beta_{f}+\beta_{g})})^2}
			=\frac{\sqrt{1+\beta_{g}\rho}-\sqrt{\rho(\beta_{f}+\beta_{g})}}
			{\sqrt{1+\beta_{g}\rho}+\sqrt{\rho(\beta_{f}+\beta_{g})}}
			$$
			and the result follows by noting that $\theta^*=\theta(\tau^*)=1-\frac{2\beta_g}{\tau^*}$.
			
			\ref{p:optirefiii}:
			Fix $\tau>0$ and $\theta\in \widehat{I}(\tau)$. We deduce from
			\cite[Example~24.2]{bauschke2011:Convex_analysis} that
			$$R_{\tau \widehat{f}}^\theta =(1+\theta)\prox_{\tau 
				\widehat{f}}-\theta \Id\colon (x_1,x_2,x_3)\mapsto 
			\left(\frac{1-\theta\tau\rho}{1+\tau\rho}x_1,
			\frac{1-\theta\tau\rho}{1+\tau\rho}x_2,
			\frac{\beta_{f}-\theta\tau}{\beta_{f}+\tau}x_3\right).$$
			Similarly, by taking $R_{\theta\tau 
				\widehat{g}}^{1/\theta}\colon (x_1,x_2,x_3)\mapsto 
			\left((\frac{\beta_{g}-\tau}
			{\beta_{g}+\theta\tau})x_1,x_2,x_3\right)$ which yields 
			$$ R_{\theta\tau 
				\widehat{g}}^{1/\theta}\circ R_{\tau \widehat{f}}^\theta \colon 
			(x_1,x_2,x_3)\mapsto 
			\left(\frac{\beta_{g}-\tau}
			{\beta_{g}+\theta\tau}\frac{1-\theta\tau\rho}{1+\tau\rho}x_1,
			\frac{1-\theta\tau\rho}{1+\tau\rho}x_2,
			\frac{\beta_{f}-\theta\tau}{\beta_{f}+\tau}x_3\right).$$
			Therefore, for every $x = (x_1,x_2,x_3) \in \HH$,
			\begin{align}
				\label{e:lastineqsharp}
				\big\|R_{\theta\tau 
					\widehat{g}}^{1/\theta}(R_{\tau 
					\widehat{f}}^\theta x)\big\|
				&\le\max\left\{\left|\frac{\beta_{g}-\tau}
				{\beta_{g}+\theta\tau}\right|\left|\frac{1-\theta\tau\rho}{1+\tau\rho}\right|
				,\left|\frac{1-\theta\tau\rho}{1+\tau\rho}\right|,
				\left|\frac{\beta_{f}-\theta\tau}{\beta_{f}+\tau}\right|\right\}
				\|x\|\nonumber\\
				&=\max\left\{\frac{1-\theta\tau\rho}{1+\tau\rho},
				\frac{\theta\tau-\beta_{f}}{\beta_{f}+\tau}\right\}
				\|x\|,
			\end{align}
			where the last equality follows from \eqref{e:leq1} and $\beta_{f} \rho \leq 1$. The equality in \eqref{e:lastineqsharp} holds by taking  $x=(0,x_2,0)$ if $\frac{1-\theta\tau\rho}{1+\tau\rho}\geq
			\frac{\theta\tau-\beta_{f}}{\beta_{f}+\tau}$ and $x=(0,0,x_3)$, otherwise.
		The case $\widehat{g}\equiv 0$ is analogous when $\beta_g=0$. The proof is complete.
	\end{proof}
	\begin{remark}
		For $\theta \in ]-1, 1+2\beta_{f}[$, set $h = f^*+\frac{1-\theta}{2} 
		\frac{\|\cdot\|^2}{2}$. By applying 
		\cite[Proposition~4.1(3)]{BricenoRoldan2026} with $\eta = 
		\frac{1-\theta}{2}$ and $\delta = 0$, we find that $h^*$ is 
		$\frac{\rho}{1+\rho\eta}$-strongly convex, $\nabla h^*$ is 
		$(\beta_{f}+\eta)$-cocoercive, and 
		$R_{\frac{1+\theta}{2}h^*}^1=R_f^\theta$. Therefore, 
		\cite[Theorem~1]{Giselsson2017IEEE} implies that $R_f^\theta$ is 
		$\rate_f(\theta)$-Lipschitz continuous, with 
		$\rate_f(\theta)$ as defined in \eqref{eq:omega_f}. 
		Note that Proposition~\ref{p:optiref}.\ref{p:optirefi} recovers this 
		same Lipschitz constant while extending its validity to any $\theta 
		\in [-1, +\infty[$.
	\end{remark}

	\begin{proposition}
		\label{t:opti}
		In the context of Problem~\ref{prob:opti}, let $\tau$ and $\theta$ in 
		$\RPP$ be such that $(1-\theta)\tau \in [-2\beta_{f},2\beta_{g}]$, let $z_0\in\HH$, and consider the recurrence
		\begin{align}
			(\forall n\in\NN)\hspace{0.2cm}
			&\left\lfloor 
			\begin{array}{ll}
				x_{n}=\prox_{\tau f}z_n\\[2mm]
				y_n=\prox_{\theta\tau g}\big((1+\theta)x_{n}
				-\theta z_n\big)\\[1mm]
				z_{n+1}=z_n+\lambda_n
				\left(\dfrac{1+\theta}{\theta}\right)
				(y_n-x_{n}),
				\label{e:RDRo}
			\end{array}
			\right.
		\end{align}
		where $(\lambda_n)_{n\in\NN}$ is a sequence in $[0,1/\alpha]$ such that 
		$\sum_{n\in\NN}\lambda_n(1-\alpha\lambda_n)=+\infty$, where
		\begin{equation}
			\alpha=
			\begin{cases}
				1,&\text{if}\:\:\beta_{f}=\beta_{g}=0;\\
				\dfrac{2\tau(1+\theta)(\beta_{f}+\beta_{g})}
				{4(\tau\beta_{g}+\beta_{f}\beta_{g}+\tau\theta\beta_{f})
					-\tau^2(\theta-1)^2},&\text{if}\:\:\beta_{f}+\beta_{g}>0,  
			\end{cases}
		\end{equation} 
		
		Then, the following hold:
		\begin{enumerate}
			\item\label{t:opti1} 
			There exists $\widehat{z}\in\Fix 
			R_{\theta\tau 
				g}^{1/\theta}R_{\tau f}^{\theta}$ such that 
			$z_n\rightharpoonup\widehat{z}$ and $x_n\rightharpoonup 
			\prox_{\tau f}\widehat{z}\in \argmin_\HH (f+g)$.
			\item\label{t:opti2}
			Suppose that $\rho>0$, $\beta_{f}+\beta_{g}>0$, and that $\lambda_n\equiv 1$.
			Then, if $(1-\theta)\tau \in ]-2\beta_{f},2\beta_{g}]$, we have
			\begin{equation}
				\label{e:linearc1}
				\|z_{n+1}-\widehat{z}\|\le \widehat{\rate}(\theta,\tau)
				\|z_n-\widehat{z}\|,
			\end{equation}
			where $\widehat{z}$ is the unique fixed point of 
			$\Fix R_{\theta\tau 
				g}^{1/\theta}R_{\tau f}^{\theta}$. Consequently,
			$(z_{n})_{n \in \NN}$ converges strongly  to $\widehat{z}$ with a linear convergence rate given by
			$\widehat{\rate}(\theta,\tau)$, where $\widehat{\rate}$ is defined in \eqref{e:r1}.
			This constant is optimized with
			$\tau^*$ and $\theta^*$ in \eqref{e:argminfg}, leading to
			$\widehat{\rate}(\theta^*,\tau^*)$ in \eqref{e:mainfg} and the rate is tight with the functions in Theorem~\ref{p:optiref}.\ref{p:optirefiii}.
		\end{enumerate}
	\end{proposition}
	\begin{proof}
		Note that \eqref{e:RDRo} is equivalent to 
		\begin{equation}
			\label{e:FPo}
			(\forall n\in\NN)\quad z_{n+1}=R_{\theta\tau g}^{1/\theta}(R_{\tau 
				f}^{\theta}z_n).
		\end{equation}
		
		\ref{t:opti1}:  This is a direct consequence of 
		Proposition~\ref{t:RDRweakc}.
		
		\ref{t:opti2}: Note that, if $(1-\theta)\tau<-2\beta_{f}$, $\widehat{\rate}(\theta,\tau)<1$ and the result follows from the Banach-Picard fixed-point theorem \cite[Theorem~1.50]{bauschke2011:Convex_analysis}. The parameters achieving the optimal rate and the functions attaining the worst case are consequence of Theorem~\ref{p:optiref}.\ref{p:optirefii}\&\ref{p:optirefiii}.
	\end{proof}
	
	\begin{remark}
		\begin{enumerate}
			\item 
			Note that, by considering relaxed resolvents, the optimal 
			convergence rate when $f$ is strongly convex and $g$ is 
			differentiable with Lipschitzian gradient is much smaller than the 
			one in \cite[Theorem~5.6]{Giselsson2017}, which is 
			$\frac{1}{1+\sqrt{\rho\beta_{g}}}$. Indeed, by setting $\beta_{f}=0$,
			$$\rate^*:=\rate(\theta^*,\tau^*)
			=\frac{1}{(\sqrt{\beta_{g}\rho+1}+\sqrt{\beta_{g}\rho})^2}<
			\frac{1}{\sqrt{\beta_{g}\rho+1}+\sqrt{\beta_{g}\rho}}
			<\frac{1}{1+\sqrt{\rho\beta_{g}}}$$
			and their values are compared in Figure~\ref{fig:comp1} in 
			Section~\ref{sec:strongmon}.
			\item 
			In the case when $\beta_{g}=0$ and $0<\rho\beta_{f}\le 1$, i.e., when 
			$f$ is strongly convex and its gradient is $\beta_{f}^{-1}$-Lipschitz,
			\eqref{e:r1} generalizes the bound in 
			\cite[Theorem~1]{Giselsson2017IEEE} allowing $\theta\neq 1$. 
			In this context, the optimal choice in our case is $\theta^*=1$ and 
			$\tau^*=\sqrt{\beta_{f}/\rho}$ recovering exactly the optimal bound 
			and parameters in 
			\cite[Proposition~3]{Giselsson2017IEEE}.
			\item The optimal Lipschitz constant in \eqref{e:mainfg} coincides with the rate of convergence of Peaceman--Rachford algorithm in \cite[Theorem~3.4(3)]{BricenoRoldan2026} without assuming the strong convexity property on $g$.
			
			\item When $\rho>0$ and $\partial f$ is $\beta_{f}$-cocoercive with $\beta_{f}>0$, the
			function $f$ is differentiable with $\beta_{f}^{-1}$-Lipschitz continuous
			gradient, and the convergence of \eqref{e:RDRo} for
			$\theta\in\left[(1+2\tau\rho)^{-1},1\right]$ and a relaxation parameter in
			$\left]0,1\right[$ also follows from
			\cite[Theorems~5.4 and~5.6]{DaoPhan2019}, where the contraction factor is
			not made explicit. Proposition~\ref{t:opti} complements these results by
			providing the constant $\rate(\theta,\tau)$ in \eqref{e:r1}, by 
			allowing
			$\theta$ up to $1+2\beta_{f}/\tau$, and by exhibiting functions for which the
			constant is attained.
			
			\item
			For merely convex $f$ and $g$, the relaxed counterpart of \eqref{e:RDRo}
			converges weakly for every $\theta\in\RPP$ whenever
			$\lambda_n\equiv\lambda<2\min\{\theta,1\}/(1+\theta)$, and this bound is
			sharp over the class of pairs in $\Gamma_0(\HH)$
			\cite[Theorem~3.3]{NilssonAkermanGiselsson2025}. Proposition~\ref{t:opti}
			addresses the complementary regime in which the cocoercivity of $\partial f$
			and $\partial g$ allows the unrelaxed recurrence \eqref{e:RDRo} and, under
			strong convexity, linear rates with explicit optimal parameters.
		\end{enumerate}
	\end{remark}
	
\begin{figure}[htbp]
		\centering
		\begin{subfigure}[b]{0.32\textwidth}
			\centering
			\includegraphics[width=\textwidth]{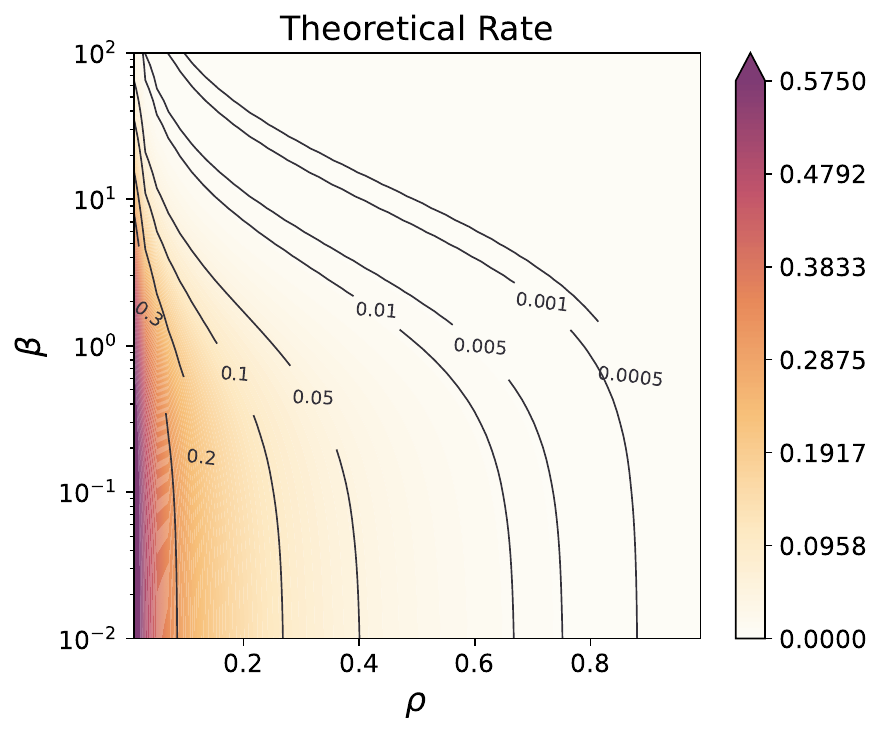}
		\end{subfigure}
		\hfill
		\begin{subfigure}[b]{0.32\textwidth}
			\centering
			\includegraphics[width=\textwidth]{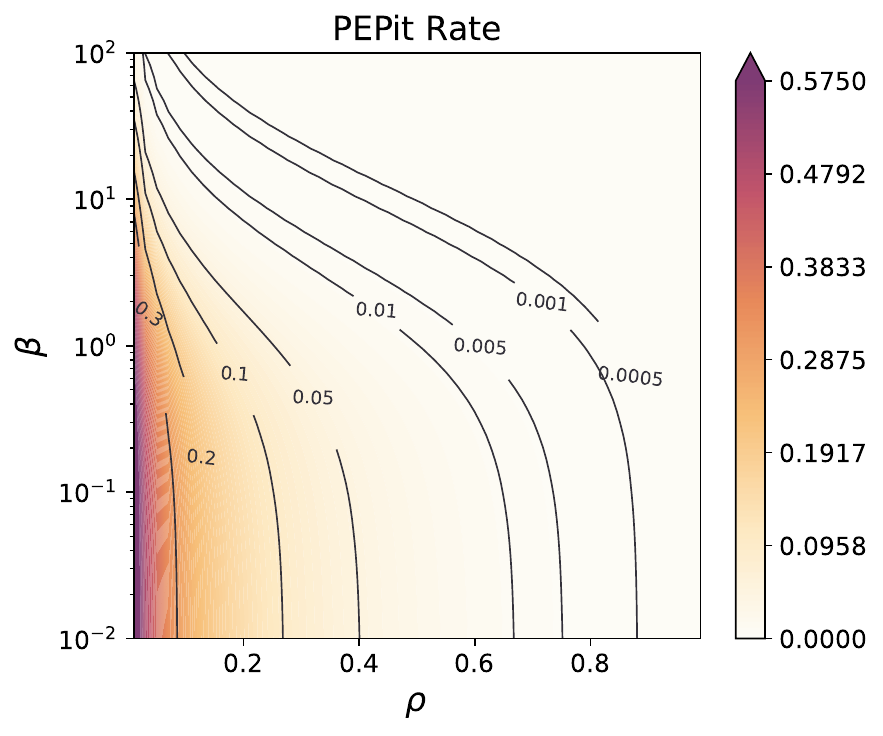}
		\end{subfigure}
		\hfill
		\begin{subfigure}[b]{0.32\textwidth}
			\centering
			\includegraphics[width=\textwidth]{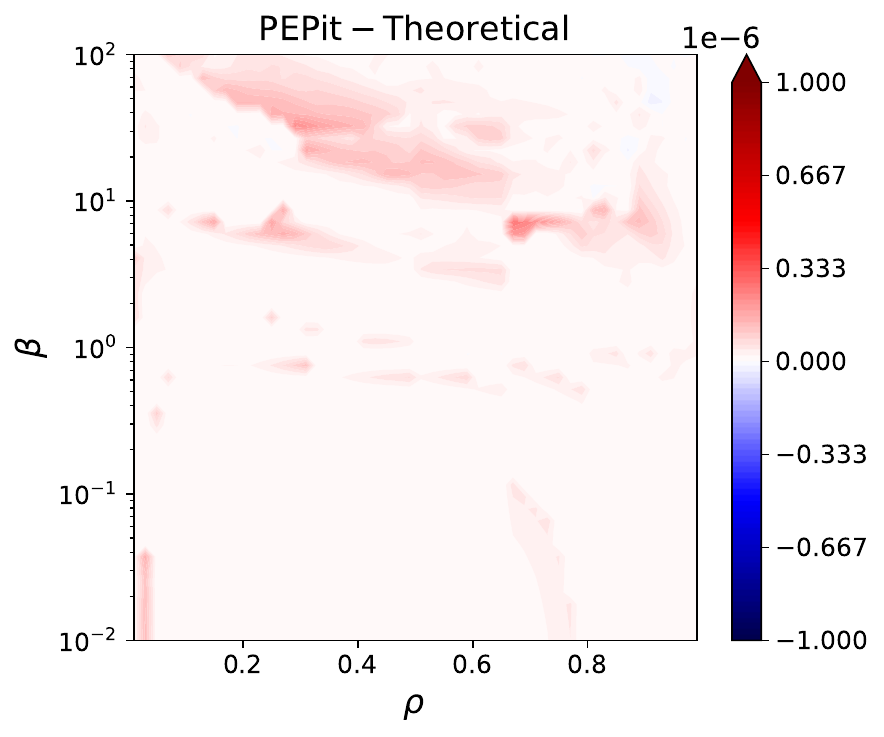}
		\end{subfigure}    
		\caption{Comparison between the theoretical convergence rate and the numerical worst-case obtained via PEPit for the Relaxed Douglas-Rachford algorithm as a function of the parameters $\rho$ and $\beta_{g}$ (log scale). We set $\beta_{f} = 1$.}
		\label{fig:rdr_comparison_convex_optim}
	\end{figure}

	We conclude the paper by comparing our theoretical rate in the convex optimization setting with 
	the computer-assisted worst-case analysis of PEPit \cite{GoujaudEtAl2024PEPit}. Figure~\ref{fig:rdr_comparison_convex_optim} shows the values of the rate obtained by our theoretical rate and that of PEPit and the difference between both values. We observe a difference of less than $10^{-6}$.

\section{Conclusion}\label{sec13}

We studied a variant of the Douglas--Rachford algorithm in which relaxed 
resolvents allow each operator to be activated with its own step-size, and we 
showed that this freedom is useful precisely when the operators are cocoercive. 
When one of them is, in addition, strongly monotone, we obtained the tight 
contraction factor of the method, its optimal parameters in closed form, and 
operators attaining the optimal rate, both for monotone inclusions and for convex 
minimization. As a by-product, we derived the tight Lipschitz constant of the 
relaxed resolvent of a strongly monotone and cocoercive operator. The resulting 
rates recover those of the Peaceman--Rachford algorithm when the second operator 
is merely monotone and improve the best known rates of classical splitting 
methods whenever it is cocoercive. All theoretical rates are corroborated by the 
numerical worst cases obtained via performance estimation.

\backmatter





\bmhead{Acknowledgements}

The work of Luis M. Brice\~no-Arias is supported by  Centro de Modelamiento Matem\'atico (CMM), FB210005, BASAL fund for 
centers of excellence,  and FONDECYT 1230257 from ANID-Chile. The work of Fernando Roldán is supported by ANID-Chile under grant FONDECYT Iniciación 11250164. Adrien Taylor is supported by the European Union (ERC grant CASPER 101162889). Views and opinions expressed are those of the authors
only and do not necessarily reflect those of the funding agencies or granting authorities, which cannot be held
responsible for them. 

\section*{Declarations}

AI tools were used to produce some figures and to numerically check tight bounds. The authors assume responsibility for all content.

\begin{itemize}
\item Funding: The work of Luis M. Brice\~no-Arias is supported by  Centro de Modelamiento Matem\'atico (CMM), FB210005, BASAL fund for 
centers of excellence,  and FONDECYT 1230257 from ANID-Chile. The work of Fernando Roldán is supported by ANID-Chile under grant FONDECYT Iniciación 11250164. Adrien Taylor is supported by the European Union (ERC grant CASPER 101162889). Views and opinions expressed are those of the authors
only and do not necessarily reflect those of the funding agencies or granting authorities, which cannot be held
responsible for them. 

\item Conflict of interest/Competing interests: The authors have no relevant financial or non-financial
interests to disclose.
\item Ethics approval and consent to participate: Not applicable
\item Consent for publication: Not applicable
\item Data and Materials availability : Not applicable
\item Code availability:  The code and experimental datasets generated or analyzed during the current study are available from the corresponding author upon reasonable request.
\item Author contribution: All the authors contributed equally to this work
\end{itemize}

\noindent

\end{document}